# An intrinsic homotopy theory for simplicial complexes, with applications to image analysis (*)

MARCO GRANDIS

*Dipartimento di Matematica, Università di Genova, Via Dodecaneso 35, 16146-Genova, Italy.*

*e-mail: grandis@dima.unige.it*

**Abstract.** A simplicial complex is a set equipped with a down-closed family of distinguished finite subsets. This structure, usually viewed as codifying a triangulated space, is used here directly, to describe "spaces" whose geometric realisation can be misleading. An intrinsic homotopy theory, not based on such realisation but agreeing with it, is introduced.

The applications developed here are aimed at image analysis in metric spaces and have connections with digital topology and mathematical morphology. A metric space $X$ has a structure $t_\varepsilon X$ of simplicial complex *at each resolution* $\varepsilon > 0$; the resulting homotopy group $\pi_n^\varepsilon(X)$ detects those singularities which can be captured by an n-dimensional grid, with edges bound by $\varepsilon$; this works equally well for continuous or discrete regions of euclidean spaces. Its computation is based on direct, intrinsic methods.

**Mathematics Subject Classifications (2000).** 55U10, 55Q05, 55U35, 18G55, 18D05, 54E35, 68U10, 54G99, 05C38.

**Key words:** Simplicial complex, homotopy groups, homotopy theory, abstract homotopy theory, 2-categories, fibre sequence, metric spaces, image processing, digital topology, digital plane, mathematical morphology, dilation, graph theory.

## Introduction

A *simplicial* (or *combinatorial*) *complex*, also called here a *combinatorial space*, is a set $X$ equipped with a family of finite subsets, the *linked parts*, such that the empty subset and all singletons are linked, and every subset of a linked part is linked. The linked parts are meant to express a notion of "proximity" or "attachment"; we shall generally avoid their classical name of *simplices*, as associated with a geometric realisation which is often inadequate for the present applications.

A *path*, or 1-*dimensional net*, in $X$ is based on a finite sequence of points $a_0,... a_n$ such that each consecutive pair $\{a_{i-1}, a_i\}$ is linked. In fact, it is a slightly more complex notion (see below), so that the set $PX$ of paths in $X$ can be equipped with a canonical combinatorial structure, yielding the path endofunctor $P$ of the category of simplicial complexes. Its powers $P^n$ allow us to define n-*tuple paths*, or n-*dimensional nets*, and n-*tuple homotopies* as maps $X \to P^nY$.

The (intrinsic) homotopy "groups" $\pi_n(X)$ of a *pointed* simplicial complex are defined, and proved to be isomorphic to the homotopy groups of the geometric realisation of $X$ (1.9; thm. 6.6). In particular, $\pi_0(X)$ is the quotient $|X|/\sim$ of the underlying pointed set, modulo the equivalence relation generated by the relation: x!x' iff {x, x'} is linked; the fundamental group $\pi_1(X)$ is the quotient of the set of loops at the base point, modulo homotopy with fixed end points.

(*) Work supported by MURST Research Projects.

The main application which we develop here is concerned with metric spaces and image analysis. If $X$ is a metric space, each real number $\varepsilon > 0$ defines a combinatorial structure $t_\varepsilon X$ on the same set, a finite part being linked iff its diameter is $\leq \varepsilon$; the space $X$ acquires thus a family of homotopy theories $\pi_n^\varepsilon(X) = \pi_n(t_\varepsilon X)$ *at resolution* $\varepsilon$, each of them detecting those singularities that can be captured by an n-dimensional net, with meshes bounded by $\varepsilon$; note that such nets are grids of points, and their moves, represented by nets of dimension $n+1$, are similarly discrete. Therefore, all this works equally well for continuous regions of $\mathbf{R}^n$ or discrete ones; in the latter case, our results are closely related with analyses of 0- or 1-connection in "digital topology" (cf. [KKM1-2]).

Consider for instance the real plane, with the *product* metric $d(\mathbf{x}, \mathbf{y}) = \vee_i \, d_i(x_i, y_i)$ (or $l_\infty$-metric, as motivated in 1.7-8), and its closed region $X$, in figure (a)

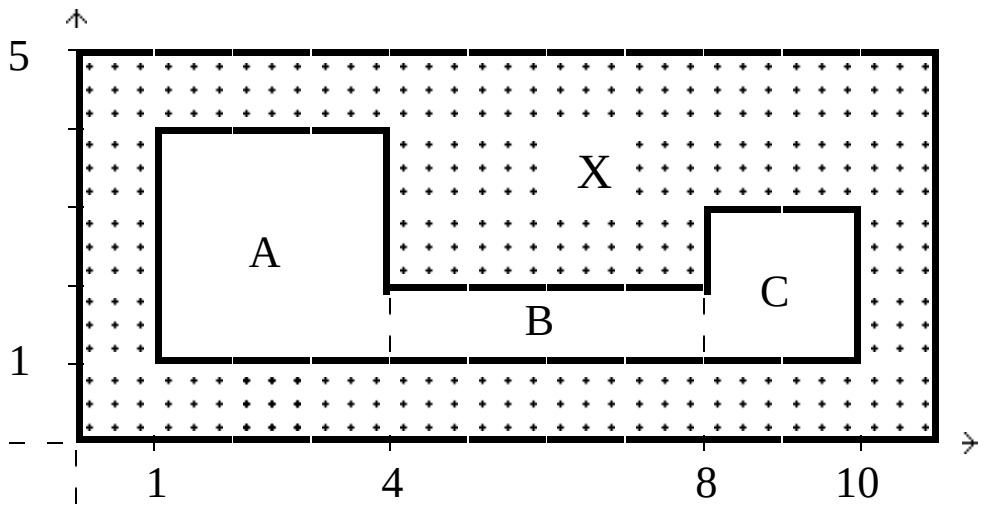

fig. (a)

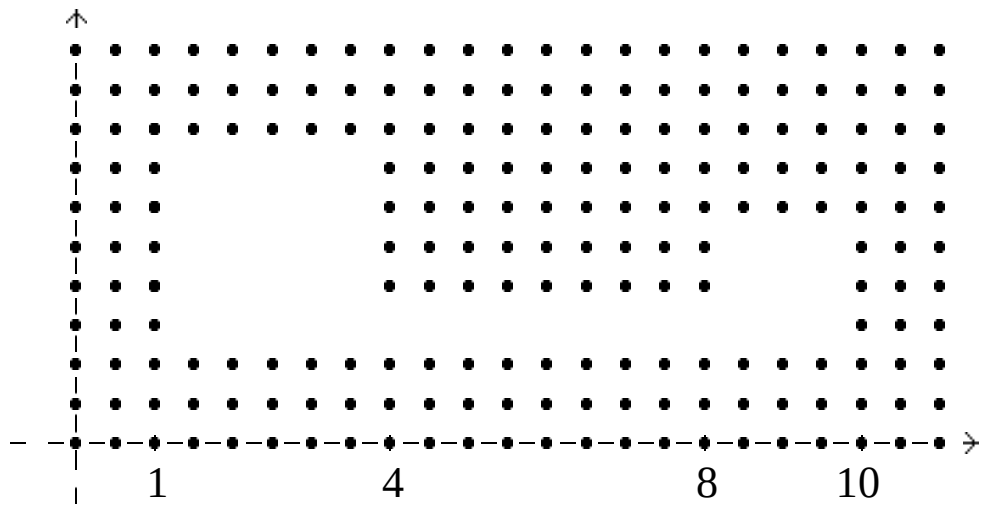

fig. (b)

Examining its fundamental group $\pi_1^\varepsilon(X)$ at increasing $\varepsilon > 0$, we distinguish: one single basin $A \cup B \cup C$ (for $0 < \varepsilon < 1$); or two basins $A$, $C$ connected by a channel $(1 \leq \varepsilon < 2)$; or one basin $A$ with a negligible appendix $(2 \leq \varepsilon < 3)$; or no relevant basin at all $(\varepsilon \geq 3)$. As shown in figure (b), this analysis can be given much in the same way on a *finite model* $X' = X \cap (\rho\mathbf{Z} \times \rho\mathbf{Z})$, as one can get from a scanning procedure at a fixed "small" resolution (e.g., $\rho = {}^1/_2$): in fact, for $\varepsilon \geq \rho$ we find the same results (while for $\varepsilon < \rho$ we have a totally disconnected object). A more detailed, but still heuristic, discussion of all this can be found in 1.8. Precise computations of the fundamental group $\pi_1^\varepsilon(X)$, based on a "van Kampen theorem" for groupoids (6.4) and a study of retracts of the combinatorial spaces $\mathbf{R}^n$ and $\mathbf{Z}^n$ (3.4-7), are given in Section 7. Higher dimensional aspects will be studied in a sequel, by higher homotopy or homology groups. Note also that, *for this* $X$, $\pi_1^\varepsilon(X)$ coincides with the topological $\pi_1$ of the *closed-spot dilation* $D_\varepsilon(X)$, the union of all closed $l_\infty$*-discs* (squares) centred at the points of $X$, with radius ${}^\varepsilon/_2$; the same holds for $X'$. This fact, for which we only have partial theoretical results (7.4-5), provides a link with "mathematical morphology" ($D_\varepsilon$ being a particular "dilation operator", cf. [He]) and with the "size homotopy groups" introduced in [FM].

A second possible application, not developed here, is concerned with a homotopy theory of "global actions", introduced by A. Bak [Ba1-2] and related to algebraic K-theory.

Let us give now some further information on simplicial complexes and their homotopy structure. Obviously, a *map* $f: X \to Y$ of such spaces is a mapping which preserves the linked parts. The resulting category **Cs** has good, well-known formal properties, which are of great help; namely, it has all limits and colimits; and it is cartesian closed ([EK]; or here, 1.4): the set of maps **Cs**$(X, Y)$ has a natural combinatorial structure, written $\mathrm{Hom}(X, Y)$ or $Y^X$, satisfying the usual exponential law, $\mathbf{Cs}(X{\times}Y, Z) = \mathbf{Cs}(X, Z^Y)$. The category $\mathbf{Cs}^*$ of *pointed simplicial complexes* (and pointed maps) also has all limits and colimits; moreover, it has a *zero-object* $0 = \{*\}$, initial and terminal.

The standard topological interval is (partially) surrogated by the *integral line* $\mathbf{Z}$ (1.3), the set of integers with the combinatorial structure of contiguity, generated by the sets $\{i, i + 1\}$ (i.e., the $t_1$-structure for $|x – y|$). A path in the simplicial complex $X$ is precisely a map $a: \mathbf{Z} \to X$ which is eventually constant at the left and at the right. The set of paths $PX \subset X^{\mathbf{Z}}$ inherits the structure of simplicial complex; the path endofunctor and its powers have a system of natural transformations comparable to the usual one, for topological spaces (Section 2), and essentially produced by the fact that $\mathbf{Z}$ is an involutive lattice *in* $\mathbf{Cs}$; but the concatenation of paths is defined by a non-natural procedure (2.5).

Homotopies are defined by maps $\alpha: X \to PY$ (3.1). They form a rather defective structure, since *consecutive homotopies cannot be concatenated*, except in particular cases, e.g. *bounded homotopies* (Section 5). The homotopy pullback $P(f, g)$ of two maps with the same codomain is easily constructed (4.5.1), but –since we do not have a general concatenation of homotopies– its 2-dimensional property can only be expressed in a restricted way. The same happens in $\mathbf{Cs}^*$; thus, the fibre sequence of a pointed map is easily obtained, by means of homotopy kernels and the loop endofunctor $\Omega$, but we lack the usual tools to study its higher dimensional properties.

We solve this problem by embedding $\mathbf{Cs}$ in a wider structure, Ps$\mathbf{Cs}$, consisting of simplicial complexes, *pseudo-maps* and *pseudo-homotopies*. Now, concatenation is possible: we obtain a lax h4-category (5.7), according to an abstract setting for homotopy developed in [Gr2] and briefly reviewed here (Section 4; 5.1). The crucial fact is that the original homotopy pullbacks in $\mathbf{Cs}$ are preserved by the embedding, and *satisfy the usual 2-dimensional property in the extension* (5.7). In Section 6, the interplay of maps, pseudo-maps and pseudo-homotopies allows us to study the fibre sequence of a pointed map and the resulting exact sequence of homotopy groups. Much in the same way as, for chain algebras, homotopy pullbacks of ordinary (multiplicative) homotopies are used to construct the fibre sequence of a morphism, but we need the more general morphisms and homotopies of chain complexes to study the higher order properties of the sequence [Gr3].

A metric simplicial complex $t_\varepsilon X$ is necessarily a *tolerance set*: it is equipped with a *tolerance* relation $x!x'$, reflexive and symmetric (defined by $d(x, x') \leq \varepsilon$), and the linked parts are the ones which are pairwise so related. Equivalently, one can consider a simple graph (as used in combinatorics); or also an *adjacency* relation, symmetric and *anti*-reflexive (as used in digital topology). The restriction of the present study to tolerance sets is possible (if $X$ is so, also $PX$ is so), but gives no real simplification. Further, there are relevant simplicial complexes which are not of this type, like the minimal circle $C_3$ and the simplicial spheres: $\mathbf{Cs}$ has no object surrogating the topological standard circle $\mathbf{S}^1$, but rather a system of k-*point circles* $C_k$, with canonical maps $C_{k+1} \to C_k$, for $k \geq 3$; similarly for higher dimensional spheres (1.3; 6.5).

Finally, some remarks on the relations with the classical theory of simplicial complexes, for which we refer to Spanier's text [Sp]. A simplicial complex (in the sense recalled above, not to be confused with simplicial set) is viewed as an elegant, economic way of codifying a triangulable space, which is reconstructed via the classical geometric realisation (1.9), gluing together a family of topological simplices (triangles of suitable dimension). There is an *intrinsic* homology theory, isomorphic to the homology of the realisation [Sp, HW] and an intrinsic *edge-path* groupoid ([Sp], 3.6), isomorphic to the fundamental groupoid of the realisation (and to the intrinsic fundamental groupoid constructed here, as proved in 2.10). On the other hand, there seems to be no general study of intrinsic homotopies, nor of intrinsic higher homotopy groups, as we are proposing here.

Here, simplicial complexes are treated as "spaces" in themselves (akin to bornological spaces, cf. 1.1). As such, they have intrinsic notions which may clash with their geometric realisation. For instance, products are not preserved, as shown by the chaotic objects **2** and **2**×**2**, realised – respectively – as a compact interval or a tetrahedron. More concretely, the geometric realisation is often misleading in the applications considered here (1.9): regions of $\mathbf{Z}^2$ and $\mathbf{Z}^3$, as produced by a regular scanning of a land or a solid object, are reconstructed as 3- or 7-dimensional spaces, respectively; continuous regions of the plane, as in figure (a), generate huge spaces (while their *spot dilations* in $\mathbf{R}^2$ or $\mathbf{R}^3$ give simple spaces having often the correct fundamental group; 7.4-5). This seems to motivate the need of intrinsic terms (like the one we are proposing, *combinatorial space*), not referring to any embedding in other categories nor leading to confusion with simplicial sets; of course, changes in a well-established terminology are difficult.

*Acknowledgements*. The homotopy theory which is being developed by A. Bak for "global actions" [Ba1-2] was a source of inspiration for this work. The author would also like to thank R. Brown and T. Porter for their hospitality at the "Workshop on global actions, groups and homotopy" (Bangor, December 1997), and the latter for providing some of the previous references.

*Notation*. The following notation for categories is used throughout: **Set** (sets); **Set*** (pointed sets); **Top** (topological spaces); **Top*** (pointed spaces); **Dm** (differential modules, i.e. differential graded R-modules, for a fixed commutative unitary ring R); **Da** (differential algebras, i.e. differential graded associative R-algebras, without unit assumption); **Gpd** (small groupoids); **Cs** (simplicial complexes); **Cs*** (pointed simplicial complexes). All these categories are equipped with the appropriate notion and structure of homotopies (cf. 4.1), except **Set**, **Set*** (or including them, with trivial homotopies). A homotopy between the maps $f, g: X \to Y$ is written as $\alpha: f \to g: X \to Y$; or also as $\alpha: X \Rightarrow Y$, without specifying $f, g$. The usual notation for intervals refers to *real* or *integral* intervals, according to the context.

## 1. The category of simplicial complexes

After recalling the basic, well known properties of this category (completeness, cartesian closedness,...) we briefly introduce the intrinsic homotopy theory and applications to be developed here. Relations with the classical theory of simplicial complexes are discussed in the Introduction and 1.9. Category theory is kept at an elementary level; undefined terms can be found in [Ma], [ASH], [Bo].

**1.1. Basic notions.** A *simplicial* (or *combinatorial*) *complex*, also called here a *combinatorial space* (*c-space* for short), is a set $X$ equipped with a *combinatorial structure*, i.e. a set $!X \subset \mathcal{P}_f X$ of finite subsets of $X$, called *linked parts*, which contains the empty subset and all singletons, and is down closed: if $\xi$ is linked, any $\xi' \subset \xi$ is so. A *morphism* of simplicial complexes, or *map*, or *combinatorial mapping* $f: X \to Y$ is a mapping between the underlying sets which preserves the linked sets: if $\xi$ is linked in $X$, so is $f(\xi)$ in $Y$. Such objects and maps form the category **Cs** of simplicial complexes.

If a structure $!X$ is contained in a second structure $!'X$ (on the same set), we say that the first is *finer* and the second is *coarser*. The combinatorial structures of a set $X$ form a complete sublattice of $\mathcal{P}(\mathcal{P}_f X)$. Any subset of $\mathcal{P}_f X$ *generates* a combinatorial structure, the finest containing it.

A *combinatorial subspace* $X' \subset X$ is a subset equipped with the induced, or initial, structure (the coarsest one making the inclusion a map): a part of $X'$ is linked iff it is so in $X$. More generally, a *subobject* $X' \prec X$ will be a subset equipped with any combinatorial structure making the inclusion a map; equivalently, $!X' \subset !X$ (this is the usual notion of *simplicial subcomplex*). The subobjects of $X$ form a complete lattice, isomorphic to a sublattice of $\mathcal{P}(!X)$; $\cap X_i$ (resp. $\cup X_i$) is the intersection (resp. union) of the underlying subsets, with structure $\cap !X_i$ (resp. $\cup !X_i$).

An equivalence relation $R$ in $X$ produces a *quotient* $X/R$, equipped with the final structure (the finest making the projection $X \to X/R$ a map): a subset of the quotient is linked iff it is the image of some linked part of $X$.

The forgetful functor $|-|$: **Cs** $\to$ **Set** has left and right adjoint, $D \dashv |-| \dashv C$:

- the *discrete* structure $DS$ is the finest one on the set $S$: the linked parts are the singletons and Ø;
- the *chaotic* structure $CS$ is the coarsest one on the set $S$: all finite subsets are linked.

The left adjoint $\pi_0 \dashv D$ is considered below (1.5). We are also interested in the category **Cs*** of *pointed* c-spaces, where an object $X$ is also equipped with a base point $*_X$, preserved by maps.

The finiteness condition on linked parts might be dropped, with some complications (cf. 2.2); but the extension of our homotopy theory to the category **Cs**' of such "generalised c- spaces", being still based on finite paths, would factor through the forgetful functor **Cs**' $\to$ **Cs**, and give no further information. (Generalised c-spaces extend *bornological spaces*, where distinguished subsets – called *bounded parts* – are also assumed to be closed under finite unions; various extensions of bornological spaces have been studied from the point of view of Categorical Topology; e.g. in [AHS], ch. 22; [Ne], ch. 21). Similarly, to drop the condition that all singletons be linked (as often done in Algebraic Topology) would produce some complications (cf. 1.4) and no further information. Finally, the requirement that the empty subset be linked is generally forgotten; but then, the empty set would have *two* combinatorial structures, which is not convenient.

**1.2. Limits.** The category **Cs** of simplicial complexes has all limits and colimits, preserved by the forgetful functor **Cs** $\to$ **Set**.

The limit $(A, p_i: A \to X_i)$ of a small diagram $X$: **I** $\to$ **Cs** is the limit of the underlying diagram $|X|$ of sets, equipped with the initial combinatorial structure: a finite subset $\alpha$ of $A$ is linked iff, for every index $i$, $p_i(\alpha)$ is linked in $X_i$. In particular, the structure of a product $\Pi X_i$ consists of the finite subsets of the products $\Pi \xi_i$ of linked parts. Analogously, the colimit $(A, u_i: X_i \to A)$ of $X$ is the colimit of $|X|$ in **Set**, with the final combinatorial structure: a subset $\alpha$ of $A$ with more than one element is linked iff it is the $u_i$-image of a linked part of $X_i$, for some index $i$. The coequaliser of two maps $f, g: X \to Y$ is thus the set-coequaliser $Y/R$, with the quotient structure.

The terminal object is the point $\{*\}$, the initial object is Ø. Moreover, $\{*\}$ represents the forgetful functor, $|-| =$ **Cs**$(\{*\}, -)$ and is a *separator*: two parallel maps $f, g: X \to Y$ are different iff there is some map $x: \{*\} \to X$ such that $fx \neq gx$. A *coseparator* is provided by the chaotic structure on two elements, $\mathbf{2} = C\{0, 1\}$. Using these two objects, it follows easily that the monomorphisms (resp. epimorphisms) of **Cs** coincide with the injective (resp. surjective) maps. A morphism $f: X$ $\to$ Y is a *regular monic* (i.e., an equaliser) iff it is injective and X has the initial structure for f: the linked parts are the ones whose image is linked; f is a *regular epi* (a coequaliser) iff it is surjective and Y has the terminal structure: the linked subsets are the images of the ones of X. Thus, a subobject $X' \prec X$ amounts to an equivalence class of arbitrary monics, while a subspace $X' \subset X$ is the same as a *regular subobject* (an equivalence class of regular monics); a quotient X/R amounts to an equivalence class of regular epis. Subspaces are classified by maps $\chi: X \to \mathbf{2}$. (It is easy to show that **Cs** is a regular category, not Barr-exact [Bo]; and also a "topological construct" [AHS, Ne].)

Also **Cs*** is complete and cocomplete. It is a *pointed* category: the singleton $0 = \{*\}$ is a zero-object, terminal and initial at the same time. Limits are calculated as in **Cs** (and suitably pointed); colimits (of non-empty diagrams) are quotients of the non-pointed corresponding ones; for instance the sum $X \vee Y$ (called a *join* in the case of pointed topological spaces) is obtained from the disjoint union by identifying the two base points; it can be realised as a subspace of the product $X \times Y$, namely $(X \times \{*\}) \cup (\{*\} \times Y)$.

**1.3. Line and spheres.** The set of integers **Z**, equipped *with the combinatorial structure of contiguity*, generated by all contiguous pairs {i, i+1}, will be called the *standard* (integral) *line* and play a crucial role in our homotopy theory. An integral interval has the induced structure, unless otherwise stated.

The structure of the *standard* (integral) n-*space* $\mathbf{Z}^n$ is generated by the "elementary cubes" $\Pi_k \{i_k, i_k+1\}$. It is relevant that the join and meet operations $\vee, \wedge: \mathbf{Z}^2 \to \mathbf{Z}$ are combinatorial mappings, as well as $-: \mathbf{Z} \to \mathbf{Z}$ (**Z** is an involutive lattice *in* **Cs**). The same holds for all translations $i \mapsto i + i_0$, whereas sum and product are *not* maps $\mathbf{Z}^2 \to \mathbf{Z}$.

The *standard elementary interval* $\mathbf{2} = [0, 1] \subset \mathbf{Z}$ is the chaotic simplicial complex on two points, C{0, 1}. The *standard elementary cube* $\mathbf{2}^n \subset \mathbf{Z}^n$ is also chaotic, as well as the *standard elementary simplex* $\mathbf{e}^n = C\{e_0,... e_n\} \subset \mathbf{Z}^{n+1}$, consisting of the unit points of the axes.

The discrete $\mathbf{S}^0 = \{-1, 1\} \subset \mathbf{Z}$ will be called the *standard 0-sphere* (and pointed at 1, when viewed in **Cs***). There is no standard circle (cf. 6.5). But, for every integer $k \geq 3$, there is a k-*point combinatorial circle*, the quotient $C_k = \mathbf{Z}/\equiv_k = \{[0], [1],... [k-1]\}$, with respect to congruence modulo k; the structure is generated by the contiguous pairs {[i], [i+1]}; the base point is [0]. (The c-spaces similarly obtained for k = 1, 2 are chaotic, hence contractible.) Such circles are related by the following maps, identifying two points

(1) $p_k: C_{k+1} \to C_k, \qquad p_k([i]) = [i] \qquad (i = 0, ... k).$

More generally, there is no standard n-sphere for $n > 0$. The *simplicial* (or tetrahedral) n-*sphere* $\Delta S^n \prec \mathbf{Z}^{n+2}$ has the same n+2 points of $\mathbf{e}^{n+1} = C\{e_0,... e_{n+1}\} \subset \mathbf{Z}^{n+2}$, but a subset is linked iff it is not total; the base point is $e_0$. The *cubical* n-sphere $\square S^n \prec \mathbf{Z}^{n+1}$ has the same $2^{n+1}$ points of the cube $\mathbf{2}^{n+1} = C\{0, 1\}^{n+1} \subset \mathbf{Z}^{n+1}$, but the linked parts are the sets of vertices contained in some face of the cube, i.e. in some hyperplane $t_i = 0$ or $t_i = 1$; the base point is 0. The *octahedral* n-sphere $\lozenge S^n = \{\pm e_0,... \pm e_n\} \subset \mathbf{Z}^{n+1}$ has 2n+2 points and the subspace structure: a subset is linked iff it does not contain opposed pairs $\pm e_i$; the base point is $e_0$. Thus, $\Delta S^0 \cong \square S^0 \cong \lozenge S^0 = \mathbf{S}^0$, $\Delta S^1 \cong C_3$, $\square S^1 \cong \lozenge S^1 \cong C_4$.

The homology of simplicial complexes is a well known tool, intrinsically defined ([Sp], ch. 4; [HW], ch. 2); let us only note that $\Delta S^n$, $\Box S^n$ and $\Diamond S^n$ are *homological n-spheres* (have the homology of the topological n-sphere), and all $C_k$ are homological circles (1.9).

**1.4. The internal hom.** The category **Cs** is cartesian closed ([EK], IV.7). The internal hom-functor $\mathrm{Hom}(X, Y) = Y^X$ is given by the set of morphisms

(1) $\mathrm{Hom}\colon \mathbf{Cs}^{op} \times \mathbf{Cs} \to \mathbf{Cs}$, $\qquad |\mathrm{Hom}(X, Y)| = \mathbf{Cs}(X, Y)$,

a finite set of maps $\varphi \subset \mathbf{Cs}(X, Y)$ being linked whenever

(2) $\varphi(\xi) = \cup_{f\in\varphi} f(\xi)$ is linked in Y, for all $\xi$ linked in X.

The usual, obvious action of Hom on morphisms does yield a combinatorial mapping

(3) $\mathrm{Hom}(u, v)\colon \mathrm{Hom}(X, Y) \to \mathrm{Hom}(X', Y')$ $\qquad (u\colon X' \to X, \ \ v\colon Y \to Y')$,

$\mathrm{Hom}(u, v)(f) = vfu$,

$(\mathrm{Hom}(u, v)(\varphi))(\xi') = v(\varphi(u(\xi'))) \in\ !Y'$ $\qquad$ (for $\xi' \in\ !X'$).

The adjunction $-\times Y \dashv \mathrm{Hom}(Y, .)$ comes from the natural bijective correspondence

(4) $\mathbf{Cs}(X\times Y, Z) \rightleftarrows \mathbf{Cs}(X, \mathrm{Hom}(Y, Z))$, $\qquad f \rightleftarrows g$,

described by the equation $f(x, y) = g(x)(y)$. Since, for $\xi\in\ !X$ and $\eta\in\ !Y$, the subsets $\xi\times\eta$ generate the combinatorial structure of $X\times Y$, it is sufficient to note that $f(\xi\times\eta) = g(\xi)(\eta)$ and that a singleton $\{x\}$ is necessarily linked. (To drop this condition would lead to a slightly more complicated internal hom and a different forgetful functor in **Set**.)

It is easy to check (but actually follows from the general theory of monoidal closed categories [Ke]) that the adjunction respects the hom-structures, providing an isomorphism (the *exponential law*)

(5) $\mathrm{Hom}(X\times Y, Z) \rightleftarrows \mathrm{Hom}(X, \mathrm{Hom}(Y, Z))$, $\qquad Z^{(X\times Y)} = (Z^Y)^X$.

Moreover, $X = \mathrm{Hom}(\{*\}, X)$. The functor $\mathrm{Hom}(Y, .)$, as a right adjoint, preserves subobjects and regular subobjects: if $Z \subset Z'$, then $\mathrm{Hom}(Y, Z) \subset \mathrm{Hom}(Y, Z')$.

(As pointed sets, **Cs*** is monoidal closed, with respect to the obvious $\mathrm{Hom}(X, Y) = \mathbf{Cs}^*(X, Y)$, pointed at the zero-map, and the *smash* product $X\wedge Y = (X\times Y)/(X\vee Y)$; this fact is not used here.)

**1.5. Connected components.** A simplicial complex X has a *tolerance relation* x!y (reflexive and symmetric), defined by $\{x, y\} \in\ !X$. The discrete functor $D\colon \mathbf{Set} \to \mathbf{Cs}$ has a left adjoint

(1) $\pi_0\colon \mathbf{Cs} \to \mathbf{Set}$, $\qquad \pi_0(X) = |X|/\!\sim$ $\qquad (\pi_0 \dashv D)$,

where $\sim$ is the equivalence relation spanned by x!y. This relation will be called *connectivity*, since a path in X (parametrised over **Z**) will be based on a finite sequence $x_1\ !\ x_2\ !\ ...\ !\ x_k$ (2.2). A non-empty c-space X is said to be *connected* if $\pi_0 X$ is a point; $\pi_0$ is called the functor of *connected components*, or *path-components* (c-spaces have no distinction for these notions, cf. 1.9). Any object is the sum of its connected components. The line **Z** is connected, as well as all its intervals.

The tolerance relation of a product $\Pi X_i$ is the product relation, $(x_i)!(y_i)$ iff $x_i!y_i$, for all i; therefore, *provided the product is finite*, the same fact holds for the connectivity relation $(x_i) \sim (y_i)$, and $\pi_0$ *preserves finite products* (but not the infinite ones, in general).

**1.6. Tolerance sets.** A *tolerance set* will be a set A equipped with a tolerance relation, i.e. a binary relation x!y in A, which is reflexive and symmetric. (Equivalently, one can assign an *adjacency* relation, *anti*-reflexive and symmetric, as more often used in digital topology [KKM1-2]; or a simple reflexive graph with vertices in A, more familiar in combinatorics.)

The obvious category **Tol** of tolerance set has mappings which preserve the tolerance relation; it is again complete and cocomplete (and cartesian closed: f!g means that x!x' implies (fx)!(gx')). Implicitly, we have already introduced (in 1.5) the forgetful functor t: **Cs** → **Tol** taking the simplicial complex X to the tolerance set tX over the same set, with x!y iff $\{x, y\} \in$ !X.

Its left and right adjoint, $d \dashv t \dashv c$, are sections of t. For a tolerance set A, dA is the finest simplicial complex on A inducing the relation ! (the non-empty linked parts are of type {x, y}, for x!y), whereas cA is the coarsest such (*a finite subset is linked iff all its pairs are* !*-related*).

We shall always *identify a tolerance set* A *with the simplicial complex* cA (*not with* dA). Thus, **Tol** becomes a full reflective subcategory of **Cs**, consisting of the c-spaces where a finite subset is linked iff all its pairs are so. The embedding c preserves all limits and is closed under subobjects; in particular, a product of tolerance sets, in **Cs**, is a tolerance set. The standard line **Z** (1.3) is a tolerance set, with i!j whenever i, j are equal or contiguous; all its powers and subobjects of powers are tolerance sets.

The embedding is not closed under quotients; in fact, the combinatorial circles $C_k = \mathbf{Z}/{\equiv_k}$ are tolerance sets for $k > 3$, but $C_3$ is not so: this simplicial complex consists of three points, linked by pairs, while the associated tolerance c-space $T = ct(C_3)$ is chaotic. (The geometric realisation takes T to a triangle, and $C_3$ to its boundary; 1.9).

More generally, the *n-truncation* $tr_n$: **Cs** → $\mathbf{Cs}_n$ $(n \geq 0)$, with values in the full subcategory of c-spaces whose linked parts have at most n+1 points, has similar left and right adjoints, *n-skeleton* and *n-coskeleton*, $sk_n \dashv tr_n \dashv cosk_n$. The cases n = 0, 1 correspond to $\mathbf{Cs}_0 \cong$ **Set**, $\mathbf{Cs}_1 \cong$ **Tol**.

**1.7. Metric spaces.** A metric space X has a family of canonical combinatorial structures $t_\varepsilon X$, *at resolution* $\varepsilon \in [0, \infty]$, where a finite subset $\xi$ is linked iff its diameter is $\leq \varepsilon$. Each of them is a tolerance set, defined by x!x' iff $d(x, x') \leq \varepsilon$. The category **Mtr** of *metric spaces and weak contractions* has thus a family of forgetful functors $t_\varepsilon$: **Mtr** → **Cs**, trivial for $\varepsilon = 0$ (giving the discrete structure) and $\varepsilon = \infty$ (the chaotic structure). The pointed case is similar. (The tolerance structures $t_\varepsilon^- X$ defined by $d(x, x') < \varepsilon$ give a less fine homotopical information; cf. 7.3).

Beware of the fact that, in **Mtr**, *a* (finite) *product has the* $l_\infty$*-metric*, given by the least upper bound $d(\mathbf{x}, \mathbf{y}) = \vee_i d_i(x_i, y_i)$; this precise metric has to be used if we want to "assess" a map, with values in a product, by its components. The functors $t_\varepsilon$: **Mtr** → **Cs** preserve finite limits; actually all limits, if metrics are allowed to take values in $[0, \infty]$, as preferable. Then $t_\varepsilon$ has a left adjoint $m_\varepsilon$: if $x \neq x'$, take $d(x, x') = \varepsilon$ for x!x', and $\infty$ otherwise; it restricts to a full embedding $m_\varepsilon$: **Tol** → **Mtr** (for $0 < \varepsilon < \infty$).

Unless differently stated, the real line **R** will have the standard metric and the tolerance structure $t_1\mathbf{R}$, with x!x' iff $|x - x'| \leq 1$, consistent with the one of **Z**; the real n-space $\mathbf{R}^n$ has the product structure, defined by the $l_\infty$-metric $\vee_i |x_i - y_i|$ (not euclidean, for $n > 1$); its linked parts are the finite subsets of all elementary cubes $\Pi_i [x_i, x_i{+}1]$; the induced structure on $\mathbf{Z}^n$ is the standard one and $\Diamond S^n = t_1\mathbf{S}^n \cap \mathbf{Z}^{n+1}$, as c-spaces. Extending what happens for the integral line (1.3), the operations of involutive lattice

(1) $\vee, \wedge: \mathbf{R}^2 \to \mathbf{R}$, $\qquad\qquad -: \mathbf{R} \to \mathbf{R}$

are weak contractions for the $l_\infty$-metric, and combinatorial mappings for all structures $t_\varepsilon$; the same holds for translations, but not for sum and product.

Marginally (at the end of Section 7), we shall also consider on the integral or real n-space $E^n$ the tolerance structure $t_\varepsilon(E^n, d_p)$ associated to the $l_p$-metric $d_p(\mathbf{x}, \mathbf{y}) = (\Sigma\, |x_i - y_i|^p)^{1/p}$ for $1 \leq p < \infty$, which is finer than the product structure $t_\varepsilon E^n = (t_\varepsilon E)^n = t_\varepsilon(E^n, d_\infty)$. For $\varepsilon = 1$, all the new structures coincide in the integral case: each point $x$ of $\mathbf{Z}^n$ is linked to the $2n$ points $x \pm e_i$, and it suffices to consider the $l_1$-metric $\Sigma\, |x_i - y_i|$. The structure $t_1(\mathbf{Z}^2, d_1)$ of the integral plane is quite different from the standard one: the diagonal $\Delta$, isomorphic to $\mathbf{Z}$ in the standard structure, becomes totally disconnected as $t_1(\Delta, d_1)$; the subset $\{0, 1\}^2$, which in the standard structure is the chaotic c-space $\mathbf{2}^2$, becomes isomorphic to the four-point circle $C_4$; $\mathbf{Z}^2$ itself is contractible, while $t_1(\mathbf{Z}^2, d_1)$ has a free fundamental group of countable rank (7.4). In digital topology, both these tolerance structures of the integral plane are used, under the names of 8-*adjacency* for the standard structure (where any point is linked to 8 others) and of 4-*adjacency* for $t_1(\mathbf{Z}^2, d_1)$.

**1.8. Homotopy at a given resolution.** As anticipated above (1.5), a path in the metric simplicial complex $t_\varepsilon X$ is based on a finite sequence of points $x_0, \ldots x_k$ with $d(x_{i-1}, x_i) \leq \varepsilon$ $(i = 1, \ldots k)$. The homotopy theory of simplicial complexes developed in the sequel will produce, for each $\varepsilon \geq 0$, a homotopy theory for pointed metric spaces, *at resolution* $\varepsilon$, $\pi_n^\varepsilon(X) = \pi_n(t_\varepsilon X)$, for which precise computations will be given in Sections 7. We sketch here some results which motivate its interest.

Loosely speaking, these groups are *concerned with singularities which can be captured by an n-dimensional net with meshes of edge* $\varepsilon$; but note that such nets are grids of points, and their moves, represented by nets of dimension $n+1$, are similarly "discrete". Of course, if $X$ is $\varepsilon$-path connected, the $\varepsilon$-homotopy groups do not depend on the base point, up to isomorphism (but the $\varepsilon$-homotopy functors do, and for them the base points cannot be ignored).

Thus, the $l_\infty$-metric space $X = T \setminus Y$, with $Y = A \cup B \cup C$ (considered in the Introduction) is $\varepsilon$-path connected as soon as $\varepsilon > 0$, but its fundamental group at resolution $\varepsilon$ varies with the latter

(1)

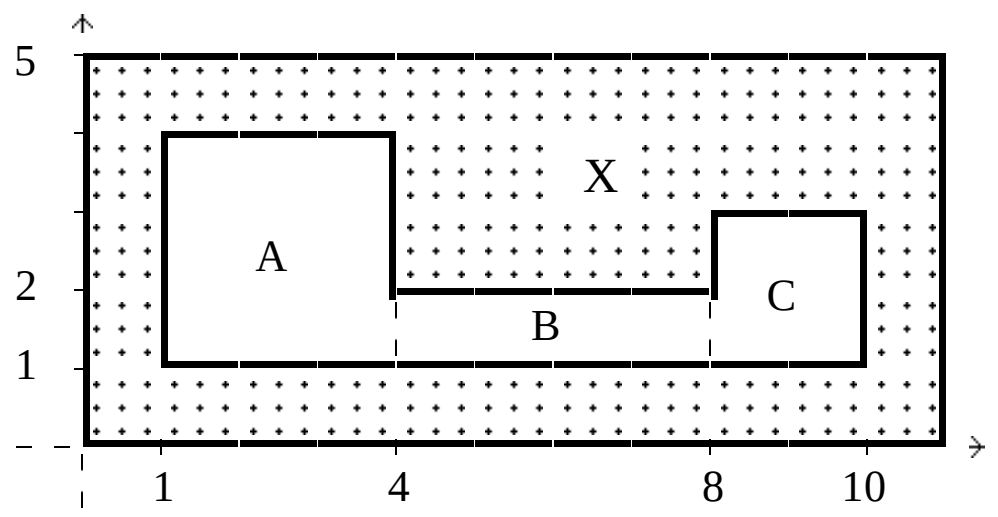

$T = [0, 11] \times [0, 5]$

$A = \,]1, 4[\, \times \,]1, 4[$
$B = [4, 8] \times \,]1, 2[$
$C = \,]8, 10[\, \times \,]1, 3[$

(2) $\pi_1^\varepsilon(X) \cong \mathbf{Z}$ $\quad (0 < \varepsilon < 1;\ 2 \leq \varepsilon < 3)$, $\qquad \mathbf{Z} * \mathbf{Z}$ $\quad (1 \leq \varepsilon < 2)$, $\qquad \{*\}$ $\quad (3 \leq \varepsilon \leq \infty)$,

detecting: one hole (corresponding to $Y$) at resolution $0 < \varepsilon < 1$; two holes (corresponding to $A$ and $C$) for $1 \leq \varepsilon < 2$ (when $B$ can be jumped over by paths); one hole again (corresponding to $A$) at resolution $2 \leq \varepsilon < 3$; and a simply connected object at resolution $\varepsilon \geq 3$. We can thus distinguish among: one single basin (or island, etc.) Y; or two basins A, C connected by a bridgeable channel (or two islands connected by an isthmus, etc.); or one basin A with a negligible appendix; or no relevant basin at all. Of course, the choice of the resolution(s) of interest should be dictated by the application (e.g., what threshold we want to fix for a lake or an island); but note that the finest description has been obtained at an intermediate resolution $(1 \leq \varepsilon < 2)$.

It is also of interest for computer graphics and image processing that our analysis of the object (1) can be given much in the same way on a finite *digital model*, as one can get from a scanning procedure at a fixed resolution $\rho$ small with respect to the dimensions of our object. Take, for instance, the trace of X on a lattice $L_\rho = \rho\mathbf{Z}\times\rho\mathbf{Z} = \{(\rho i, \rho j) \mid i, j \in \mathbf{Z}\}$ at resolution $\rho = k^{-1}$, for an integer $k \geq 2$ (as in figure (b) of the Introduction). The metric space $X' = X \cap L_\rho$ is totally disconnected at resolution $\varepsilon < \rho$; it is $\varepsilon$-path connected for $\varepsilon \geq \rho$, where the group $\pi_1^\varepsilon$ gives the same results as above

(3) $\pi_1^\varepsilon(X') \cong \mathbf{Z}$ $(\rho \leq \varepsilon < 1;\ 2 \leq \varepsilon < 3)$, $\mathbf{Z} * \mathbf{Z}$ $(1 \leq \varepsilon < 2)$, $\{*\}$ $(3 \leq \varepsilon \leq \infty)$.

For $\rho = 1$, all this is still true, but the first case is empty. (If $\rho$ is not of type $k^{-1}$, these results have a "marginal" variation, due to the interference of $L_\rho$ with the boundary of X in $\mathbf{R}^2$. This effect is rather artificial, due to a *hybrid* definition of X' as the trace of a given "continuous space" on $L_\rho$. Naturally, as well as practically, we should rather start from an explicit description of X' in terms of points of the lattice, as one would get from a real scanning.)

Similarly, the left hand (resp. the central) figure below, a metric subspace of $\mathbf{R}^2$, is viewed by the fundamental $\varepsilon$-group as a circle (resp. a "figure 8") at resolution $1 \leq \varepsilon < 8$, then as a trivial object

(4)

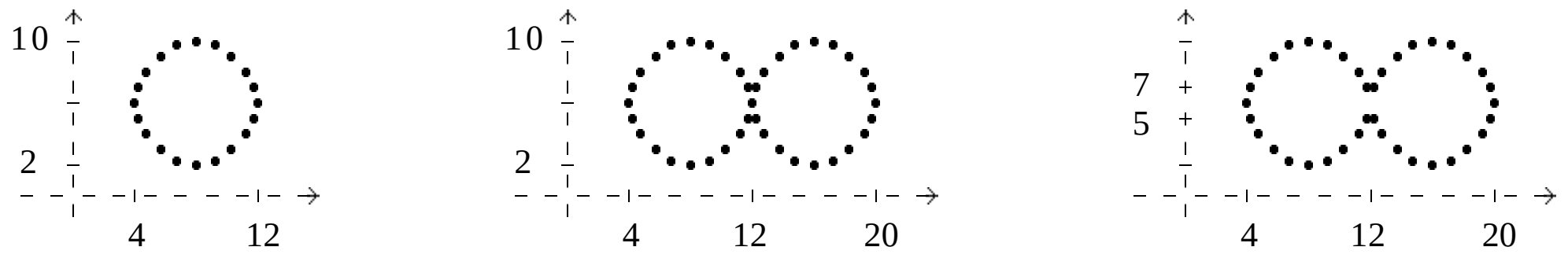

while the right hand figure is analysed as a circle for $1 \leq \varepsilon < 2$, and a "figure 8" for $2 \leq \varepsilon < 8$; a small $\varepsilon$ is sensitive to "errors".

Higher dimensional aspects should be studied by higher homotopy or homology groups; this is deferred to a sequel. Finally, let us note that the euclidean metric of the plane, being invariant by rotation, might seem to be more adequate for the present applications. This is not necessarily true, since a scanning procedure can introduce privileged directions, as above. Computation is easier in the $l_\infty$-metric, where one can take full advantage of cartesian products (cf. 3.4-7; 7.1-2), yet also possible in the $l_p$-metrics (7.4-5).

**1.9. Topological realisations.** A c-space X has a *simplicial* geometric realisation $\mathfrak{S}(X)$ ([Sp], ch. 3), obtained by gluing together a topological simplex (triangle) $\Delta^n(\xi)$ of dimension n, for each linked part $\xi$ having n+1 elements, so to preserve the inclusion of linked parts (cf. 6.6).

(a) This is adequate to the present purposes for various "elementary" simplicial complexes: for instance, $\mathfrak{S}(\mathbf{Z})$ is the gluing of a sequence of edges, homeomorphic to the real line; $\mathfrak{S}(C_k)$ is the boundary of a k-gon, homeomorphic to $\mathbf{S}^1$; the simplicial sphere $\Delta S^n$ (1.3) gives the boundary of the topological simplex $\Delta^{n+1}$, which is homeomorphic to $\mathbf{S}^n$.

(b) The geometric realisation is complicated and geometrically inadequate for c-spaces *of cubical type* (loosely speaking), like $\mathbf{Z}^n$ or the cubical sphere $\square S^n$, as soon as $n \geq 2$. Each elementary square of $\mathbf{Z}^2$ (a chaotic c-space on four points) is turned into a (solid) tetrahedron; the integral plane is realised as a pasting of tetrahedra along edges, a sort of "bubble wrap" formed of (solid, deformed) tetrahedra

(1) 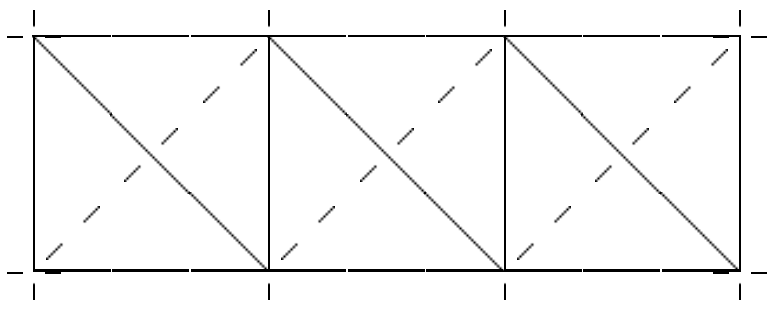

Similarly, the geometric realisation of $\square S^2$ is not a 2-sphere, but six tetrahedra pasted along the edges of a cube; $\mathbf{Z}^n$ is turned into a topological space of dimension $2^n - 1$. Homotopically, all this is already a problem, as the realisation may be of scarce help in computing the homotopy groups. Geometrically, it is even worse; regions of $\mathbf{Z}^2$ and $\mathbf{Z}^3$, as produced by a scanning of a land or a solid object, are important for image analysis; we do not want to reconstruct them as 3- or 7-dimensional spaces, respectively.

(c) Less "tame" c-spaces, like $t_1\mathbf{R}$ or $t_1\mathbf{S}^1$, even having a quite simple homotopy type, may have huge geometric realisations.

The much simpler "spot dilations" of a metric subspace $X \subset \mathbf{R}^n$ (7.4-5), notwithstanding their limited range, appear to be of more interest for the simplicial complexes considered in (b), (c).

Finally, let us also recall the "finite type" representation $X \mapsto \mathcal{A}(!X)$ introduced by McCord ([Mc]; see also [HaV] and references therein). Start from the embedding $!: \mathbf{Cs} \to \mathbf{Pos}$, which associates to a simplicial complex $X$ the (pre)ordered set $!X$ of its linked parts and compose it with the embedding $\mathcal{A}: \mathbf{Pos} \to \mathbf{Top}$ which endows a preordered set $S$ with the *Alexandroff* topology: opens are the down-closed subsets $U$ of $S$ ($x \leq x' \in U$ implies $x \in U$); $\mathcal{A}$ is an isomorphism onto the category of *Alexandroff spaces* (topological spaces where open sets are stable under intersection). Simplicial complexes are thus embedded in the category of locally finite $T_0$ Alexandroff spaces, and *finite objects are turned into finite spaces*. For instance, $\mathcal{A}(!C_k)$ has $2k$ points, corresponding to the vertices and edges of a k-gon; the vertices are open, while the least neighbourhood of an edge consists of itself and its two vertices. The homotopy groups of $\mathcal{S}(X)$ and $\mathcal{A}(!X)$ coincide [Mc]; they also coincide with the intrinsic homotopy groups of $X$ introduced here (6.6). Note also that, in agreement with the terminology of 1.5, all spaces $\mathcal{S}(X)$ and $\mathcal{A}(!X)$ are locally path-connected, whence connected *iff* path-connected.

## 2. Paths in simplicial complexes

We study now the paths of **Cs**; the fundamental groupoid derived from such path (2.9) is proved to be isomorphic to the edge-path groupoid (2.10). The letter $\kappa$ denotes a Boolean variable, with values $\pm 1$ (often contracted to $\pm$).

**2.1. Lines.** The standard line $\mathbf{Z}$ has been equipped with the combinatorial structure of contiguity, generated by the contiguous pairs $\{i, i+1\}$ (1.3). A map $a: \mathbf{Z} \to X$ amounts to a sequence of points in the c-space $X$, written $a(i)$ or $a_i$, with $a_i \,!\, a_{i+1}$ for all $i \in \mathbf{Z}$. The *simplicial complex of lines* of X is $L(X) = Hom(\mathbf{Z}, X) = X^{\mathbf{Z}}$; a finite set $\Lambda$ of lines is linked iff each set $\cup_{a\in\Lambda} \{a_i, a_{i+1}\}$ is linked in X (for $i\in\mathbf{Z}$). Note that the *set* |LX| only depends on the tolerance relation in X, but this is not true of the object LX (unless X itself is a tolerance set).

We have a representable endofunctor $L = Hom(\mathbf{Z}, –)$: **Cs** → **Cs**. The iterated functor $L^2$: **Cs** → **Cs** is represented as $Hom(\mathbf{Z}^2, –)$, by the *plane* $\mathbf{Z}^2 = \mathbf{Z}\times\mathbf{Z}$, whose linked sets are generated by the "elementary squares" $\{i, i+1\}\times\{j, j+1\}$, and induced by the product tolerance relation on $\mathbf{Z}^2$: (i, j) ! (i', j') iff (i!i' and j!j'). $L^2(X)$ is the *c-space of planes* of X. The *set* $|L^n(X)|$ depends on the linked 2n-tuples of X.

The standard line and its powers are connected by relevant maps: *degeneracy* (e), *connections* $(g^-, g^+)$, *reversion* (r), and *interchange* (s)

(1) $\{*\} \xleftarrow{\ e\ } \mathbf{Z} \overset{g^\kappa}{\leftleftarrows} \mathbf{Z}^2$ $\qquad$ $r: \mathbf{Z} \to \mathbf{Z}$, $\qquad$ $s: \mathbf{Z}^2 \to \mathbf{Z}^2$

$g^-(i, j) = i\vee j$, $\qquad$ $g^+(i, j) = i\wedge j$,

$r(i) = -i$, $\qquad$ $s(i, j) = (j, i)$,

giving **Z** the structure of an involutive lattice in **Cs** (as already observed in Section 1). Applying the contravariant functor Hom(–, X), we get the corresponding natural transformations: degeneracy $e: 1 \to L$, connections $g^\kappa: L \to L^2$, reversion $r: L \to L$ and interchange $s: L^2 \to L^2$.

**2.2. Paths.** A *path*, or 1-*dimensional net*, of the simplicial complex X is a line $a: \mathbf{Z} \to X$ eventually constant at the left and at the right: there is a finite interval $\rho = [\rho^-, \rho^+] \subset \mathbf{Z}$ $(\rho^- \leq \rho^+)$ such that a is constant on the half-lines $]-\infty, \rho^-]$, $[\rho^+, \infty[$

(1) $a(i) = a(\rho^-)$ for $i \leq \rho^-$, $\qquad$ $a(i) = a(\rho^+)$ for $i \geq \rho^+$,

and determined by its values over $\rho$; the latter will be called an (admissible) *support* of a and often viewed as an element of $J \subset |\mathbf{Z}|^2$, the *set* of increasing pairs of integral numbers, with no combinatorial structure. Note that a path is not assigned any precise support (this would prejudice the deriving homotopies, cf. 4.6). An *immediate*, or *one-step* path has support **2** = [0, 1] and amounts to a pair of linked points $(a_0, a_1) \in X^2$.

The end points, or *faces*, of the path a, $\partial^\kappa a = a(\rho^\kappa)$, are well defined. The *path object* $PX \subset X^{\mathbf{Z}}$ is the combinatorial subspace of paths. The *path functor* P: **Cs** → **Cs** acts on a morphism $f: X \to Y$ as a subfunctor of $(–)^{\mathbf{Z}}$

(2) $Pf: PX \to PY$, $\qquad$ $(Pf)(a) = fa$,

taking a to a path which admits the same supports, and possibly smaller ones.

The faces are *maps* $\partial^\kappa: PX \to X$ (and natural transformations $P \to 1$): in fact, each (*finite*) linked set of paths admits a common support $\rho$, and $\partial^\kappa$ is the evaluation at the point $\rho^\kappa$ for all such paths. (Dropping the finiteness condition on linked parts would require here a non trivial modification; the simplest one is perhaps to replace **Z** with the *extended* integral line $\overline{\mathbf{Z}} = \mathbf{Z} \cup \{\pm\infty\}$.)

Two points $x, x' \in X$ are linked by a path in X iff $x \sim x'$, for the equivalence relation generated by the tolerance relation of X (1.5). The quotient yields the functor of connected components, $\pi_0$: **Cs** → **Set**, $\pi_0 X = X/\sim$, left adjoint to the discrete functor D (1.5).

In **Cs***, the path functor is the previous one, pointed at the constant path $0_{*_X}$ (of the base point $*_X$ of X); the functor of path-components takes values in pointed sets, $\pi_0$: **Cs*** $\to$ **Set***. Note that, as for pointed topological spaces, the paths of X are *not* pointed maps, nor has the line **Z** been pointed (we are still using the cartesian closed structure of **Cs**).

**2.3. Double paths.** A *square*, or *double path*, or 2-*dimensional net,* in X is an element of $P^2(X)$.

The inclusion $PX \subset LX$ (preserved by L, a right adjoint) gives $P^2X \subset LP(X) \subset L^2(X)$. Thus, a square amounts to a plane $A: \mathbf{Z}^2 \to X$ admitting a *support* $\rho\times\sigma = [\rho^-, \rho^+] \times [\sigma^-, \sigma^+] \subset \mathbf{Z}^2$

(1) $A(i, j) = A(\rho^\kappa, j)$, for $\kappa i \geq \kappa\rho^\kappa$,

$A(i, j) = A(i, \sigma^\kappa)$, for $\kappa j \geq \kappa\sigma^\kappa$,

which means that $A: \mathbf{Z}^2 \to X$ is constant on the dotted lines and shaded quadrants below

(2) 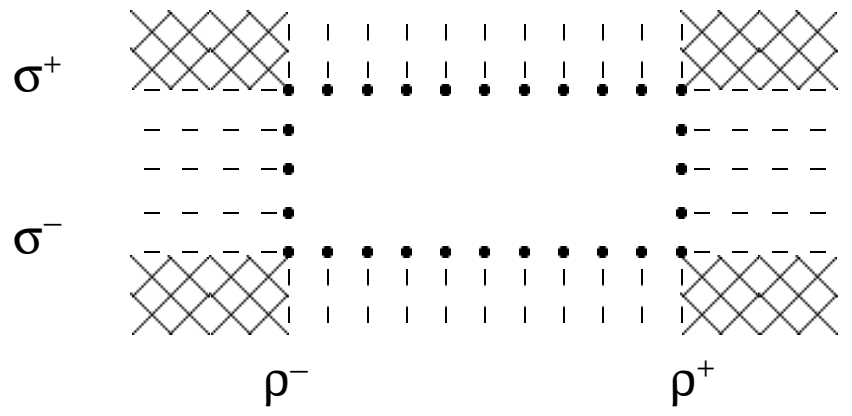

and determined by its values over $\rho\times\sigma$. The two faces $\partial^\kappa: P \to 1$ of the path functor produce four faces of $P^2$, $\partial^\kappa P: P^2 \to P$ and $P\partial^\kappa: P^2 \to P$, which we view as *horizontal* and *vertical*, respectively. Thus, the faces of $A \in P^2X$ are four paths, which connect four vertices

(3) $(\partial^\kappa P.A)(i) = A(i, \sigma^\kappa)$, $(P\partial^\kappa.A)(j) = A(\rho^\kappa, j)$,

$$\begin{array}{ccc} y' & \xrightarrow{\;b\;} & z \\ {\scriptstyle u}\uparrow & A & \uparrow{\scriptstyle v} \\ x & \xrightarrow[\;a\;]{} & y \end{array} \qquad \begin{array}{ll} a = \partial^- P.A, & b = \partial^+ P.A \\ u = P\partial^-.A, & v = P\partial^+.A \\ x = \partial^-.\partial^- P.A = \partial^-.P\partial^-.A, \ ... & \end{array}$$

Double paths whose vertical faces are trivial are of special interest (2.8). The object of n-*dimensional paths* $P^nX$ consists of the maps $a: \mathbf{Z}^n \to X$ admitting a *support* $\rho_1\times...\times\rho_n \subset \mathbf{Z}^n$; it has 2n faces, $\partial_i^\kappa = P^{n-i}\partial^\kappa P^{i-1}: P^n \to P^{n-1}$ (i = 1,... n; $\kappa = \pm$).

**2.4. The structure of paths.** Since P is a subfunctor of $L = (-)^{\mathbf{Z}}$ and $P^2$ of $L^2$ (2.3), other natural transformations (after faces) are inherited from the ones of L (2.1.1), which can plainly be restricted: degeneracy (e), connections $(g^-, g^+)$ and symmetries (r, s)

(1) $$1 \ \underset{e}{\overset{\partial^\kappa}{\leftleftarrows\!\!\rightarrow}} \ P \ \overset{g^\kappa}{\rightrightarrows} \ P^2 \qquad\qquad r: P \to P, \qquad s: P^2 \to P^2;$$

in particular, every point $x \in X$ has a *trivial path* $0_x = e(x)$ (constant at x), and every path $a \in PX$ from x to x' has a *reversed path* $-a = r(a)$, from x' to x, namely $(-a)(i) = a(-i)$.

These seven natural transformations satisfy the axioms of an (associative) *cubical comonad with symmetries* [Gr1, 3] (the functorial dual-analogue of a *commutative, involutive cubical monoid,* i.e. a set equipped with two structures of commutative monoid $g^\kappa$, where the unit $(\partial^\kappa)$ of each operation is an absorbent element for the other, and the involution r turns each structure into the other):

(2) $\partial^\kappa.e = 1$, $g^\kappa.e = Pe.e$ $(= eP.e)$ (*degeneracy axiom*),

$Pg^\kappa.g^\kappa = g^\kappa P.g^\kappa$, $P\partial^\kappa.g^\kappa = 1 = \partial^\kappa P.g^\kappa$ (*associativity, unit*),

$P\partial^\kappa.g^\eta = e.\partial^\kappa = \partial^\kappa P.g^\eta$ (*absorbency*; $\kappa \neq \eta$),

$r.r = 1$, $r.e = e$, $\partial^-.r = \partial^+$, $g^-.r = Pr.rP.g^+$,

$s.s = 1$, $s.Pe = eP$, $P\partial^\kappa.s = \partial^\kappa P$, $s.g^\kappa = g^\kappa$,

$Pr.s = s.rP$ (*symmetries*).

It will be useful to note that P preserves all *finite* limits of simplicial complexes, but does not preserve infinite products. In fact, $(-)^{\mathbf{Z}}$ preserves all limits (as a right adjoint); now, a *finite* jointly monic family of maps $X \to X_i$ (e.g., the projections of a limit) reflects the lines which admit a support; this obviously fails in the infinite case. Therefore, P has no left adjoint, and the homotopies which it generates have no cylinder functor.

The same holds in $\mathbf{Cs}^*$, where one shows in the same way that the cocone functor K and the loop functor $\Omega$ (consisting of the loops at the base-point)

(3) $KX = \mathrm{Ker}\,(\partial^-: PX \to X) = \{a \in PX \mid \partial^-a = *_X\}$,

$\Omega X = \mathrm{Ker}\,(\partial^-: PX \to X) \cap \mathrm{Ker}\,(\partial^+: PX \to X) = \{a \in PX \mid \partial^-a = *_X = \partial^+a\}$,

preserve *finite* limits and do not have a left adjoint (there is no cone nor suspension). The interchange $s: P^2 \to P^2$ restricts to isomorphisms

(4) $\overline{s}: \Omega P \to P\Omega$, $s_\Omega: \Omega^2 \to \Omega^2$, $(\partial^\kappa\Omega.\overline{s} = \Omega\partial^\kappa: \Omega P \to \Omega)$,

as is obvious from the concrete description of $\Omega$, but can also be formally deduced from the previous preservation property (cf. 4.5.2 for a similar argument).

**2.5. Path concatenation.** X is always a simplicial complex. We want now to concatenate two *consecutive* paths a, b $(\partial^+a = \partial^-b)$; this cannot be done in a natural way. Any choice of a pair $\rho(a)$, $\rho(b)$ of admissible supports for a, b yields an *admissible concatenation* c of our paths, with "pasting point" at $\rho^+(a) + \rho^-(b)$

(1) $c(i) = a(i - \rho^-(b))$, for $i \leq \rho^+(a) + \rho^-(b)$,

$c(i) = b(i - \rho^+(a))$, for $i \geq \rho^+(a) + \rho^-(b)$,

and admissible support $\rho(a) + \rho(b) = [\rho^-(a) + \rho^-(b), \rho^+(a) + \rho^+(b)]$. To verify that c is indeed a map, take a linked part $\zeta = \{i, i+1\} \subset \mathbf{Z}$; then, we can compute $c(\zeta)$ via a (if $i+1 \leq \rho^+(a) + \rho^-(b)$) or b (if $i \geq \rho^+(a) + \rho^-(b)$), up to a translation in $\mathbf{Z}$ (an automorphism of $\mathbf{Cs}$). We show below that *two admissible concatenations* c, c' *of* a, b *are always congruent up to delays*, $c \equiv c'$ (2.6), to deduce later that they are *homotopic with fixed end points* (2.8).

In order to get a well defined operation, from now on $\rho(a)$ will denote the *standard support* of the path a, i.e. the least admissible one (for a constant path, we always take [0, 0]); the (standard) *concatenation* a + b will be the one produced by this selection. Thus, $\rho(0_x) = 0$, $\rho(-a) = [-\rho^+(a), -$ $\rho^-(a)]$ and $\rho(a+b) = \rho(a) + \rho(b)$. The *set* |PX|, with concatenation and reversion, is a (true!) category with involution (over |X|, and in additive notation)

(2) $(a + b) + c = a + (b + c)$, $\quad 0 + a = a = a + 0$,

$- 0 = 0$, $\quad - (a + b) = - b - a$, $\quad - (- a) = a$.

A map $f: X \to Y$ takes a path $a \in PX$ to a path $fa \in PY$ admitting $\rho(a)$ as *a support*; thus $f.(a + b)$, computed by means of $\rho(a)$, $\rho(b)$, is an admissible concatenation of fa, fb and

(3) $f.(a + b) \equiv fa + fb$.

Note that our operation is based on a mapping *of sets* $\rho: |PX| \to J$, which is *not natural* (for maps $f: X \to Y$); this double anomaly is further discussed in 2.11.

**2.6. Congruence up to delays.** We want to show how the various admissible concatenations of two paths (2.5) are related. The *monoid* $D_1$ *of delays* is defined as the submonoid of $\mathbf{Cs}(\mathbf{Z}, \mathbf{Z})$ generated by the family of *elementary delays* $\delta_t$ $(t \in \mathbf{Z})$

(1) $\delta_t(i) = i$, if $i \leq t$, $\quad \delta_t(i) = i - 1$, otherwise.

The crucial point is a sort of *simplicial identity*

(2) $\delta_t.\delta_{s+1} = \delta_s.\delta_t$ $\quad (t \leq s)$,

and the following consequence, *the main property of delays*, or *cofiltering propery*

(3) for any two delays $d_1, d_2$, there are delays $e_1, e_2$ such that $d_1e_1 = d_2e_2$.

(This holds *within* the generators by the simplicial identity; the general case follows by a repeated application of the particular one.) A delay is not a path in **Z**; but if $a: \mathbf{Z} \to X$ is a path, so is any delayed line $ad: \mathbf{Z} \to X$. Two paths $a, b \in PX$ are said to be *congruent* (up to delays), $a \equiv b$, if there exist two delays $d, d': \mathbf{Z} \to \mathbf{Z}$ such that $ad = bd'$. Congruence is an equivalence relation, by the main property above. Note that *any path is congruent to its translations*: the path $a(i-1)$ can be obtained as $a\delta_t(i)$, for any $t \leq \rho^-(a)$. (On the other hand, if $t \geq \rho^+(a)$, then $a\delta_t = a$.)

Now, it is easy to show that *all admissible concatenations of two given paths* a, b *are congruent*: let c be one of them, derived from supports $\rho, \sigma$; then, varying $\rho^-$ or $\sigma^+$ has no effect on c, while increasing $\rho^+$ (resp. $\sigma^-$) of one unit yields a concatenation $c' = c\delta_t$ delayed of one unit at a suitable instant. One can easily prove that $\equiv$ is consistent with concatenation, but we shall not need this.

In dimension 2, we use similarly the submonoid $D_2$ of $\mathbf{Cs}(\mathbf{Z}^2, \mathbf{Z}^2)$ generated by the elementary delays $\delta_t \times \mathbf{Z}$, $\mathbf{Z} \times \delta_{t'}$ $(t, t' \in \mathbf{Z})$. Since the generators of the first system commute with the ones of the second, the main property of delays still holds, and we have a 2-*congruence* relation $A \equiv_2 B$ in $P^2X$.

**2.7. Pasting.** For a double path $A: X \to P^2Y = P(PY)$ (2.3), our selection of supports already determines *one* standard support $(\rho(A), \sigma(A)) \in J^2$ (2.2), *consistent with interchange*: in fact A (resp. sA) is a 1-dimensional path *in* PY, between its faces $\partial^\kappa P.A$ (resp. $P\partial^\kappa.A$).

Two *vertically consecutive* squares $A, B: X \to P^2Y$ $(\partial^+P.A = \partial^-P.B)$ have a *vertical* (standard) *concatenation*, or *vertical pasting* $C = A +_v B$ (as consecutive paths between their horizontal faces $X \to PY$), which depends on the "vertical" components $\sigma = \sigma(A)$, $\tau = \sigma(B)$ of their standard supports

(1) $C(i, j; x) = A(i, j - \tau^-; x)$, for $j \leq \sigma^+ + \tau^-$,

$C(i, j; x) = A(i, j - \sigma^+; x)$, for $j \geq \sigma^+ + \tau^-$;

note that the lower *horizontal face* of a square is its *vertical domain*, i.e. the domain for the vertical sum. The horizontal faces of C are obvious, the vertical ones are *an admissible* concatenation of the vertical faces of A and B, *congruent* to the standard one

(2) $\partial^-P.(A +_v B) = \partial^-P.A$, $\partial^+P.(A +_v B) = \partial^+P.B$,

$P\partial^\kappa.(A +_v B) \equiv P\partial^\kappa.A + P\partial^\kappa.B$.

Symmetrically, two *horizontally consecutive* squares $A, B: X \to P^2Y$ ($P\partial^+.A = P\partial^-.B$) have a *horizontal concatenation* (or *horizontal pasting*)

(3) $A +_h B = s(sA +_v sB)$.

Both operations are associative, satisfy the cancellation property, have identities and involution

(4) $0_v(a) = eP.a$, $0_h(a) = Pe.a$ (*vertical and horizontal identity of* a),

(5) $-_v A = rP.A$ $-_h A = Pr.A$ (*vertical and horizontal reversion of* A).

Plainly, all admissible pastings of two given squares are 2-congruent (up to 2-dimensional delays, 2.6). It follows that our two operations satisfy the *four middle interchange law, up to 2-congruence*

(6) $(A +_h B) +_v (C +_h D) \equiv_2 (A +_v C) +_h (B +_h D)$;

(choose a support $(\rho, \sigma)$ admissible for all four squares; performing the above operations with this uniform choice, we would get an equality; the standard supports give 2-congruent results).

**2.8. Homotopy of paths.** More particularly, a *2-path* is a double path whose vertical faces are trivial. Such squares are stable under horizontal and vertical reversion or pasting.

A 2-path is viewed as a *homotopy with fixed end points*, or 2-*homotopy*, $A: a \simeq_2 b$, between two paths, its horizontal faces ($a = \partial^-P.A$, $b = \partial^+P.A$); plainly, a and b have the same end points x, x' and the vertical faces of A are $0_x = P\partial^-.A$, $0_{x'} = P\partial^+.A$

(1) $A(i, \sigma^-) = a(i)$, $A(i, \sigma^+) = b(i)$, $A(\rho^-, j) = x$, $A(\rho^+, j) = x'$.

(We shall speak of an *immediate*, or *one-step* 2-*homotopy*, when $A(-, j)$ coincides with a for $j \leq 0$, with b for $j \geq 1$.) The resulting relation $a \simeq_2 b$ is an equivalence, as follows easily from the vertical structure considered above:

- reflexivity: given a path a, take its vertical identity $0_v(a): a \simeq_2 a$;
- symmetry: given $A: a \simeq_2 b$, take its vertical reverse $-_vA = rP.A: b \simeq_2 a$;
- transitivity: given also $B: b \simeq_2 c$, take their vertical sum $A+_vB: a \simeq_2 c$.

A crucial fact is that *two congruent paths* $a \equiv b$ (2.6) *are always 2-homotopic*, $a \simeq_2 b$ (*Caterpillar Lemma*). The elementary case is proved by the caterpillar homotopy $a\delta_t \simeq_2 a$, constructed below; since we already know that $\simeq_2$ is an equivalence relation, the conclusion follows. Therefore, all admissible concatenations of two paths are 2-homotopic; moreover, in 2.7.2, $P\partial^\kappa.(A +_v B) \simeq_2 P\partial^\kappa.A + P\partial^\kappa.B$.

*The Caterpillar Homotopy*. Consider a path $a: \mathbf{Z} \to X$ and the delayed path $b = a\delta_t: \mathbf{Z} \to X$ (displayed below), produced by the elementary delay $\delta_t$ (2.6.1).

There is a homotopy with fixed end points $A: b \simeq_2 a$, which modifies $b$ by a sort of caterpillar progress, taking the delay at the right of the support of $a$, where it is ineffective (the "caterpillar wave" is in boldface letters)

| | | | | | | | | | | |
|---|---|---|---|---|---|---|---|---|---|---|
| | … | $t$ | $t+1$ | $t+2$ | … | $s-1$ | $s$ | $s+1$ | … | $(i \in \mathbf{Z})$ |
| a: | … | $a_t$ | $a_{t+1}$ | $a_{t+2}$ | … | $a_{s-1}$ | $\mathbf{a}_s$ | $\mathbf{a}_s$ | … | $(j = s)$ |
| | … | $a_t$ | $a_{t+1}$ | $a_{t+2}$ | … | $\mathbf{a}_{s-1}$ | $\mathbf{a}_{s-1}$ | $a_s$ | … | |
| (2) | … | | | | | | | | … | |
| | … | $a_t$ | $\mathbf{a}_{t+1}$ | $\mathbf{a}_{t+1}$ | … | $a_{s-2}$ | $a_{s-1}$ | $a_s$ | … | $(j = t+1)$ |
| b: | … | $\mathbf{a}_t$ | $\mathbf{a}_t$ | $a_{t+1}$ | … | $a_{s-2}$ | $a_{s-1}$ | $a_s$ | … | $(j = t)$ |

we form thus a finite sequence of !-linked paths, from $b$ to $a$, of length $s - t + 1$ (where $s \geq t \vee \rho^+(a)$), and globally a square $A: b \simeq_2 a$ $(A(i, j) = a\delta_{j \vee t}(i))$.

(*Note that the "one-step try" is wrong*: taking $B(i, j) = b(i)$ for $j \leq 0$, and $a(i)$ otherwise, does not yield a map $\mathbf{Z}^2 \to X$, because the subsets $\{a_i, a_{i+1}, a_{i+2}\}$ need not be linked.) One can similarly build a 2-homotopy $a(\delta_t)^r \simeq_2 a$ starting from an r-delayed path, with a caterpillar wave of length $r+1$. But this follows for free from the elementary case, for *all* delays in $D_1$, as we know that $\simeq_2$ is an equivalence relation.

Finally, let us briefly sketch the analogous result one dimension up: if the squares $A, B: \mathbf{Z}^2 \to X$ are 2-congruent (end of 2.6), then they are also *homotopic with fixed faces*: there is a triple path in $P^3X$ having two parallel faces equal to $A$, $B$, the others being degenerate. Again, it suffices to prove the elementary case, $B = A.(\delta_t \times \mathbf{Z})$ or $B = A.(\mathbf{Z} \times \delta_t)$; but this follows for free from cartesian closedness and the 1-dimensional result, viewing $A$ as a path $\mathbf{Z} \to X^{\mathbf{Z}}$, in the two possible ways.

**2.9. The fundamental groupoid.** We have now all what we need to define the (intrinsic) *fundamental groupoid* of a simplicial complex, $\Pi X = |PX| / \simeq_2$.

As a quotient, $\Pi X$ is a category over the object-set $|X|$, with respect to concatenation $[a] + [b] = [a + b]$. Moreover, $[a]$ and $[-a]$ are reciprocal; to prove that $-a + a \simeq_2 0$, it is sufficient to take the horizontal pasting of $A = \mathrm{Pr}.g^-(a)$ and $B = g^-(a)$

$$\text{(1)} \qquad \begin{array}{ccccc} x' & \xrightarrow{\;0\;} & x' & \xrightarrow{\;0\;} & x' \\ {\scriptstyle 0}\uparrow & A & \uparrow{\scriptstyle a} & B & \uparrow{\scriptstyle 0} \\ x' & \xrightarrow[-a]{} & x & \xrightarrow[a]{} & x' \end{array}$$

(2) $-a + a \;=\; \partial^-P.A + \partial^-P.B \;\simeq_2\; \partial^-P(A +_h B) \;\simeq_2\; \partial^+P(A +_h B) \;=\; 0_{x'}$.

Globally, we have the (intrinsic) *fundamental groupoid functor*, with values in the category of small groupoids

(3) $\Pi: \mathbf{Cs} \to \mathbf{Gpd}$, $\qquad \Pi X = |PX| / \simeq_2$, $\qquad (\Pi f)[a] = [fa]$.

In fact, given a 2-path A: $a \simeq_2 a'$, the 2-path fA: $\mathbf{Z}^2 \to X$ shows that $fa \simeq_2 fa'$; consistence with composition is already known, by the formula $f.(a + b) \equiv fa + fb$ (2.5.3) and the previous Caterpillar Lemma. Further properties of this functor (as homotopy invariance) have to be deferred after the study of homotopies (6.2). It can be noted that ΠX is determined by the four-tuples of linked points of X; the result below shows that, actually, linked triplets are sufficient.

**2.10. The edge-path groupoid.** A simplicial complex X has a known, intrinsic *edge-path* groupoid $\mathcal{E}X$ ([Sp], 3.6), isomorphic to the fundamental groupoid of the geometric realisation. We recall its construction, with some minor simplifications, and prove that it is isomorphic to the intrinsic fundamental groupoid ΠX, derived from the path functor P.

To begin with, an *edge path* in X is precisely a finite non empty sequence $a = (a_0,... a_h)$ of points of X, with $a_{i-1}!a_i$, and goes from $\partial^- a = a_0$ to $\partial^+ a = a_h$. Given a consecutive $b = (b_0,... b_k)$, with $a_h = b_0$, their concatenation is $a + b = (a_0,... a_h, b_1,... b_k)$. We have obtained a category EX with objects the points of X and an obvious involution (or reversion): $-a = (a_h,... a_0)$. Two edge paths a, b are *simply equivalent* if they only differ by the repetition of a point, or can be "deformed one into the other along a linked triplet", i.e. if they can be expressed as below

(1) $c + (x, y, z) + d, \qquad\qquad c + (x, z) + d,$

(or symmetrically) where c, d are paths and {x, y, z} is a linked subset. The *equivalence* relation ~ generated by simple equivalence is a congruence of involutive categories, and the quotient $\mathcal{E}X = EX/\!\sim$ is a groupoid, the *edge path groupoid* of the simplicial complex X.

Edge paths have a precise support (which makes them easy to manage), while P-paths do not (which is necessary for other reasons, discussed in 4.6). We assign to an edge path $a = (a_0,... a_h)$ its obvious extension $\hat{a} \in PX$, *admitting* the support [0, h]

(2) $\hat{a}: \mathbf{Z} \to X, \qquad\qquad \hat{a}(i) = a_{(0\vee i)\wedge h}\,.$

This mapping $(-)\hat{}: |EX| \to |PX|$ preserves faces and identities; reversion and concatenation are preserved up to delays: $(-a)\hat{} \equiv -\hat{a}$, $(a + b)\hat{} \equiv \hat{a} + \hat{b}$. If a is obtained from b by repeating a point, then $\hat{a}$ is obtained from $\hat{b}$ by a delay; if a = (x, y, z) is a linked triplet and b = (x, z), then $(x, y, z)\hat{} \simeq_2 (x, x, z)\hat{}$, by a one-step homotopy with fixed endpoints. We have thus an induced functor on groupoids, $(-)\hat{}: \mathcal{E}X \to \Pi X$, which is the identity on objects and surjective on arrows.

To verify it is an isomorphism, it suffices to consider two paths $\hat{a} \simeq_2 \hat{b}$ and prove that a ~ b. We can assume that a, b have the same "length" $h \geq 2$ and that $\hat{a} \simeq_2 \hat{b}$ is a *one-step* 2-homotopy, so that all four-tuples $\{a_i, a_{i+1}, b_i, b_{i+1}\}$ are linked. If a, b only differ at one index i, then $0 < i < h$ and a ~ b because

(3) $(a_{i-1}, a_i, a_{i+1}) = (a_{i-1}, a_i, b_{i+1}) \sim (a_{i-1}, a_i, b_i, b_{i+1}) \sim (a_{i-1}, b_i, b_{i+1}) = (b_{i-1}, b_i, b_{i+1});$

the general case follows from the previous one:

(4) $a = (b_0, a_1, a_2,... a_h) \sim (b_0, b_1, a_2,... a_h) \sim ... \sim (b_0, b_1, b_2,..., b_{h-1}, a_h) = b.$

**2.11. Remarks on concatenation.** The concatenation of consecutive paths has been defined by a selection of standard supports, a *mapping of sets* ρ: |PX| → J which is *not natural*. Its formal ("algebraic") realisation, considered below, is similarly defective.

First, the object QX *of (pairs of) consecutive paths*, or Q-*pullback* of X, is a pullback in **Cs**

(1) $QX = PX \times_X PX = \{(a, b) \in PX{\times}PX \mid \partial^+(a) = \partial^-(b)\}$,

(2) $$\begin{array}{ccc} PX & \xrightarrow{\partial^+} & X \\ {\scriptstyle k^-}\uparrow & \overset{k^+}{\searrow} & \uparrow{\scriptstyle \partial^-} \\ QX & \longrightarrow & PX \end{array} \qquad\qquad k^-(a, b) = a, \qquad k^+(a, b) = b,$$

with projections $k^\kappa$: QX → PX; it has three faces QX → X, namely $\partial^{--} = \partial^-k^-$ (*lower*), $\partial^{\pm} = \partial^+k^- = \partial^-k^+$ (*middle*), $\partial^{++} = \partial^+k^+$ (*upper*).

The concatenation of consecutive paths forms now a family of *set-mappings* kX: |QX| → |PX|, *natural up to congruence* (2.5.3)

(3) kX: |QX| → |PX|, k(a, b) = a + b,

Pf.kX ≡ kY.Qf (f(a + b) ≡ fa + fb),

which satisfies the standard properties (as considered in [Gr4]): $\partial^\kappa k = \partial^\kappa k^\kappa$, $ke_Q = e$, $kr_Q = rk$, kP.s' = s.Pk (where $e_Q$: 1 → Q, $r_Q$: Q → Q and s': PQ → QP are induced by e, r, s respectively; the last property is concerned with the pasting of double paths, 2.7). (The "normal" situation studied in [Gr4] has a *natural* family of *morphisms* kX: QX → PX, granting the concatenation of homotopies; naturality can be weakened up to $\simeq_2$, but the second anomaly is not easily overcome.)

## 3. Homotopies of simplicial complexes

Homotopies in **Cs** and **Cs*** are introduced. In the integral and real spaces $\mathbf{Z}^n$, $\mathbf{R}^n$, we construct an essential tool, telescopic homotopies and telescopic retracts.

**3.1. Homotopies.** A *homotopy* of simplicial complexes α: f → g: X → Y is a map α: X → PY such that $\partial^-\alpha = f$, $\partial^+\alpha = g$. It can also be viewed as a map α: X → $Y^{\mathbf{Z}}$, or α: **Z**×X → Y, such that every line α(x) admits a support ρ(x) ∈ J and

(1) α(i, x) = f(x), for $i \leq \rho^-(x)$, α(i, x) = g(x), for $i \geq \rho^+(x)$.

Our homotopy is said to be *bounded* if it admits a constant support ρ(x) = ρ; and *bounded on connected components* if this holds on every connected component of X. Similarly, in a *left bounded* (resp. *positive*, *immediate*) homotopy α, all paths α(x) admit a support $[\rho^-, \rho^+(x)]$ (resp. $[0, \rho^+(x)]$, [0, 1]). An immediate homotopy amounts to a map α: **2**×X → Y, or also to two linked maps f!g in $Y^X$: for each ξ∈ !X, f(ξ)∪g(ξ) is linked in Y (this relation f!g is called "contiguity" in [Sp], 3.5); if Y is a tolerance set, this reduces to f(x)!g(x'), for all x!x' in X.

The category **Cs** will always be equipped with these (general) homotopies and the operations produced by P, as a cubical comonad with symmetries (2.4.2):

(a) *whisker composition* of maps and homotopies (for u: X' → X, v: Y → Y'):

v∘α∘u: vfu → vgu (v∘α∘u = Pv.α.u: X' → PY'),

(b) *trivial homotopies*: $0_f: f \to f$ $(0_f = ef: X \to PY)$,

(c) *reversion*: $-\alpha: g \to f$ $(-\alpha = r\alpha: X \to PY)$.

The whisker composition is also written $v\alpha u$ when no ambiguity may arise. The functor P is extended to homotopies, *by means of the interchange* s, which turns the faces $P\partial^\kappa Y$ into the required faces $\partial^\kappa PY$ of P(PY)

(2) $\hat{P}(\alpha) = s.P\alpha: PX \to P^2Y$, $\partial^\kappa P.s.P\alpha = P\partial^\kappa.P\alpha = Pf^\kappa$ .

By the structural properties of P (2.4.2), this extension preserves trivial homotopies and reversion: $\hat{P}(0_f) = s.P(ef) = eP.Pf = 0_{Pf}$, $\hat{P}(-\alpha) = -\hat{P}(\alpha)$.

A *double homotopy* is a map $\Phi: X \to P^2(Y)$. Equivalently, it is a map $\Phi: \mathbf{Z}^2{\times}X \to Y$ such that every square $\Phi(-; x)$ admits a support $(\rho(x), \sigma(x)) \in J^2$ (2.3)

(3) $\Phi(i, j; x) = \Phi(\rho^\kappa, j; x)$, for $\kappa i \geq \kappa\rho^\kappa(x)$,

$\Phi(i, j; x) = \Phi(i, \sigma^\kappa; x)$, for $\kappa j \geq \kappa\sigma^\kappa(x)$.

As for double paths (2.3), its four faces are homotopies, connecting four maps

(4) $(\partial^\kappa P.\Phi)(i, x) = \Phi(i, \sigma^\kappa(x); x)$ (*horizontal faces*),

$(P\partial^\kappa.\Phi)(j, x) = \Phi(\rho^\kappa(x), j; x)$ (*vertical faces*).

A double homotopy $\Phi: X \to P^2Y = P(PY)$ is an ordinary homotopy of its horizontal faces $\partial^\kappa P.\Phi: X \to PY$, whereas the interchanged $s\Phi$ is an ordinary homotopy of the vertical faces of $\Phi$, $P\partial^\kappa.\Phi$. We say that $\Phi$ is a *bounded* double homotopy if both these ordinary homotopies are bounded, i.e. there is a constant support $(\rho, \sigma) \in J^2$.

A finite family of homotopies $\alpha_i: f_i \to g_i: X_i \to Y_i$ has a *product* $\Pi\alpha_i: \Pi f_i \to \Pi g_i$, represented by the product map, $\Pi\alpha_i: \Pi X_i \to \Pi(PY_i) = P(\Pi Y_i)$ (2.4).

All this holds similarly for pointed objects, where homotopies are defined by (pointed) maps $\alpha: X \to PY$. Now, also the loop-functor $\Omega$ is extended to homotopies, by means of the interchange $\bar{s}: \Omega P \to P\Omega$ (2.4.4)

(5) $\hat{\Omega}(\alpha) = \bar{s}.\Omega\alpha: \Omega X \to P(\Omega Y)$, $\partial^\kappa\Omega.\bar{s}.\Omega\alpha = \Omega\partial^\kappa.\Omega\alpha = \Omega f^\kappa$.

**3.2. The homotopy relation.** By the previous structure, the *homotopy relation* $f \simeq g$, defined by the existence of a homotopy $f \to g$, is a reflexive and symmetric relation, "weakly" consistent with composition ($f \simeq g$ implies $vfu \simeq vgu$), but *presumably* not transitive.

In fact, in the simplicial complex $Y^X$, we have various generalised "homotopy relations":

(1) $f!g \Rightarrow f \simeq_b g \Rightarrow f \simeq_c g \Rightarrow f \simeq g \Rightarrow f \simeq_t g \Rightarrow f \simeq_p g$.

(i) We know that $f!g$ iff there is an immediate homotopy $f \to g$ (3.1).

(ii) The equivalence relation $f \sim g$ spanned by the former amounts therefore to the existence of a *bounded homotopy* $f \to g$, and will be preferably written as $f \simeq_b g$; it is a congruence of categories (an equivalence relation, consistent with composition).

(iii) The existence of a homotopy $f \to g$ bounded on connected components, $f \simeq_c g$, is also so.

(iv) The relation $f \simeq g$ is strictly weaker. Actually, $\mathrm{id}\mathbf{Z}$ is homotopic to the constant map $0: \mathbf{Z} \to \mathbf{Z}$ (by *telescopic homotopy*, 3.4); but such maps are not $\simeq_b$-related (a path from 0 to n has at least length n). The transitive relation spanned by $\simeq$ will be written $\simeq_t$; it is a congruence in **Cs**.

(v) Finally, the *pseudo-homotopy* relation $f \simeq_p g$ will mean that, for each *finite* subset $K \subset X$, the restrictions $f_K, g_K: K \to Y$ are homotopic (whence $\simeq_b$-homotopic); it is again a congruence and will be studied in Section 5.

The functor $P: \mathbf{Cs} \to \mathbf{Cs}$ is *homotopy invariant*, in the sense that any homotopy $f \to g$ yields a homotopy $Pf \to Pg$ (3.1.2); the same holds for $P, K, \Omega: \mathbf{Cs}^* \to \mathbf{Cs}^*$.

The map $f: X \to Y$ is a *homotopy equivalence* if it has a homotopy inverse g ($gf \simeq 1$, $fg \simeq 1$). We shall say that two objects are *homotopy equivalent*, $X \simeq Y$, if they are linked by a *finite sequence* of homotopy equivalences; this implies that they are isomorphic objects modulo $\simeq_t$. The object X is *contractible* if it is homotopy equivalent to a point.

On a *finite* object X, the previous relations (ii)-(v) coincide; we generally write $f \simeq f'$ this equivalence relation, and [f] the homotopy class of a map f. A wider class of objects over which $\simeq$ is transitive is considered below (3.3).

**3.3. Positive retracts.** A *deformation retract* S of a simplicial complex X is a subspace whose inclusion u has a retraction p, with $up \simeq \mathrm{id}X$

(1) $u: S \rightleftarrows X: p,$ $\qquad pu = 1, \quad \alpha: up \to 1_X,$

and we speak of a *positive* (resp. *bounded*, *immediate*) deformation retract when the homotopy $\alpha$ can be so chosen (3.1).

Thus, an immediate deformation retract $S \subset X$ is a subspace whose inclusion u has a retraction p, with $(\mathrm{id}X)!(up)$, i.e. $\xi \cup up(\xi)$ is linked for all $\xi \in\, !X$. X is immediately contractible to its point $x_0$ iff the latter can be added to any linked part (or is linked to any point, in a tolerance set). A non-empty chaotic space is immediately contractible.

A relevant information on X is the existence of a *finite, positive deformation retract* S. Then, we shall prove that two maps $f, g: X \to Y$ are homotopic iff their restrictions $fu, gu: S \to Y$ are so (for any Y, cf. 5.2); *the homotopy relation* $f \simeq g$ *is therefore an equivalence in any set* $\mathbf{Cs}(X, Y)$ (and a congruence in the full subcategory of the objects admitting such retracts).

More particularly, X is *positively contractible* if it admits a positive deformation retract reduced to a point; or equivalently, if there is a positive homotopy $\alpha: c \to 1_X$ starting from a constant endomap. Then each pair of maps $f, g: X \to Y$ is homotopic, $f \simeq g$. The objects having a finite positive deformation retract and the positively contractible objects are closed under finite products (use a product homotopy, 3.1). All this still holds replacing *positive* with *immediate*.

The special classes of deformation retracts which we construct in the rest of this Section are essential for the theory and will allow us to calculate fundamental groupoids, in Section 6 (by van Kampen). They can also give a better understanding of combinatorial homotopies, for which the usual intuition based on topology may deceive.

**3.4. Telescopic homotopies.** First, we show that $\mathbf{Z}^n$, $\mathbf{R}^n$, and wide classes of their regions are positively contractible. The need of the following constructions will be better understood noting that the usual topological homotopy $\varphi: [0, \rho]\times\mathbf{R}^2 \to \mathbf{R}^2$, $\varphi(t, x) = t.\rho^{-1}.x$ $(\rho > 0)$ *is not a combinatorial mapping*, with respect to any structures $t_\varepsilon\mathbf{R}^2$, $t_{\varepsilon'}[0, \rho]$ $(\varepsilon, \varepsilon' > 0)$; in fact, consider the linked subset $\xi = \{x_0, x_1\} \subset \mathbf{R}^2$, where $x_0 = (1, 1)$ and $x_1 = (1 - \varepsilon, 1 + \varepsilon)$, and note that $\varphi([1, 1+\varepsilon']\times\xi)$ cannot be linked, since $x_0$ is the only point of the main diagonal whose distance from $x_1$ is $\leq \varepsilon$.

Let E denote either the integral line **Z**, or the real line **R**. There is a rather obvious *telescopic homotopy* (centred at the origin) constructed by means of the lattice operations (1.3, 1.7)

(1) $\alpha: 0 \to \mathrm{id}: E \to E$,

$$\alpha(i, x) = 0\vee(i\wedge x), \qquad \alpha(i, -x) = -\alpha(i, x) \qquad (x \geq 0),$$

and *positive*: its path $\alpha(-, x)$ has a positive (standard) support, namely $[0, \rho^+(x)]$, with $|x| \leq \rho^+(x) < |x|+1$ (the proof that $\alpha$ is indeed a map $\mathbf{Z}\times E \to E$ will be given in 3.6.2, for a more general construction). Similarly, any integral or real (or rational) interval is positively contractible; and so is any finite product of such intervals.

This homotopy can be viewed as a collection of "telescopic arms" of height $\rho^+(x)$, opening for increasing $i \geq 0$

(2)

| | | | | | | | | | |
|---|---|---|---|---|---|---|---|---|---|
| … | −3 | −2 | −1 | 0 | 1 | 2 | 3 | … | $(x\in\mathbf{Z})$ |
| … | −2 | −2 | −1 | 0 | 1 | 2 | 2 | … | $(i = 2)$ |
| … | −1 | −1 | −1 | 0 | 1 | 1 | 1 | … | $(i = 1)$ |
| … | 0 | 0 | 0 | 0 | 0 | 0 | 0 | … | $(i = 0)$ |

Interestingly, there is *no* positive homotopy in the opposite direction $\mathrm{id} \to 0$ (which suggests that this "dual notion of positive retract" is of scarce interest). In fact, in the integral case, any positive homotopy $\beta: \mathrm{id} \to f: \mathbf{Z} \to \mathbf{Z}$ is necessarily trivial, since the map $g = \beta(1, -)$ adjacent to id

(3)

| | | | | | | | |
|---|---|---|---|---|---|---|---|
| … | $g(-2)$ | $g(-1)$ | $g(0)$ | $g(1)$ | $g(2)$ | … | $(i = 1)$ |
| … | $-2$ | $-1$ | $0$ | $1$ | $2$ | … | $(i = 0)$ |

must coincide with the former ($g(j)$ is linked with $j-1, j, j+1$, whence $g(j) = j$), and so on. The same holds in the real case: examine the maps $\beta(i, -): \mathbf{R} \to \mathbf{R}$, for $0 \leq i \leq 1$. It follows that *the only* bounded *deformation retract of the integral or real line is the line itself.*

For the n-dimensional space $E^n$, a *telescopic homotopy* will be any product of 1-dimensional telescopic homotopies (centred at any point) and trivial homotopies. For instance, for $n = 2$, consider $\beta = \alpha\times\alpha$ (*centred at the origin*) and $\gamma = \alpha\times 0_{id}$ (*centred at the second axis*)

(4) $\beta: 0 \to \mathrm{id}: E^2 \to E^2$, $\qquad \beta(i, x_1, x_2) = (\alpha(i, x_1), \alpha(i, x_2))$,

(5) $\gamma: p_2 \to \mathrm{id}: E^2 \to E^2$, $\qquad \gamma(i, x_1, x_2) = (\alpha(i, x_1), x_2)$.

The first can be visualised displaying, at the left hand, the motion of the point $x = (4, 2)$

(6) 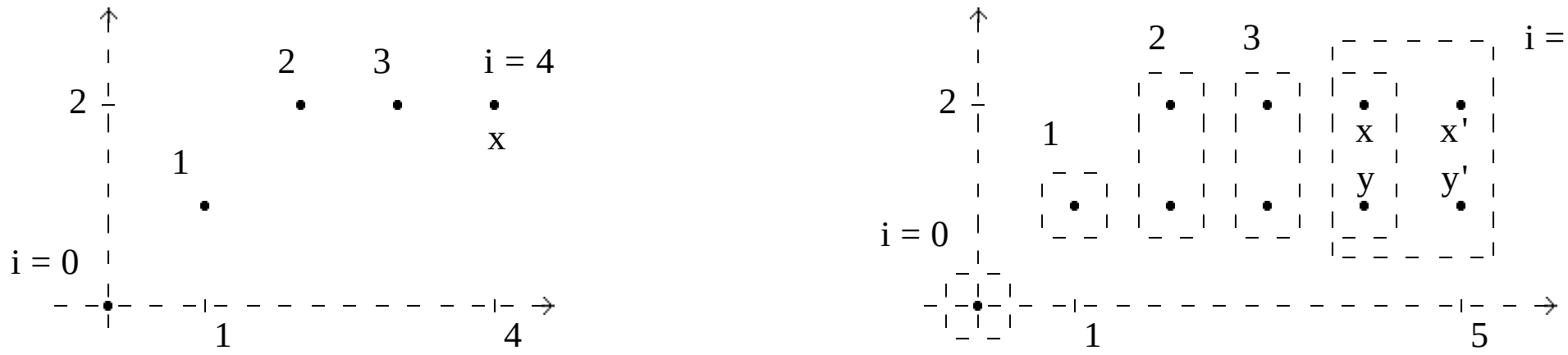

along $\beta$, i.e. the path $\beta(-, x)$: $0 \to y$, from the instant $i = 0$, when $x$ is sent to $0(x) = 0$, to the instant $i = 4 = \max(4, 2)$ when $id(x) = x$ is reached; the path follows the diagonal of the first quadrant up to $i = 2 = \min(4, 2)$. At the right hand, we display similarly the progress of a linked part $\xi = \{x, x', y, y'\}$, from the instant $i = 0$ (when $0(\xi) = \{0\}$), to the instant $i = 5$ when $id(\xi) = \xi$ is reached; it is easy to verify, directly, that all subsets $\beta(\{i, i+1\}{\times}\xi)$ are linked.

**3.5. Telescopic retracts.** Telescopic homotopies can be restricted to regions of the integral or real spaces. Let us illustrate some useful cases in the plane $E^2$, either integral ($\mathbf{Z}^2$) or real ($\mathbf{R}^2$).

Let us start from the telescopic homotopy $\beta = \alpha{\times}\alpha$: $0 \to id$, centred at the origin (3.4.4). Say that a subspace $X \subset E^2$ is *telescopically contractible to the origin* if, for any $x \in X$, the path $\beta(-, x)$, from the origin to $x$, is contained in $X$. Then, $X$ admits the origin as a positive deformation retract, via the restriction of $\beta$. This notion is analogous to a star shaped region, for topological homotopy; but a telescopically contractible domain may look quite different, as in the left hand example below

(1) 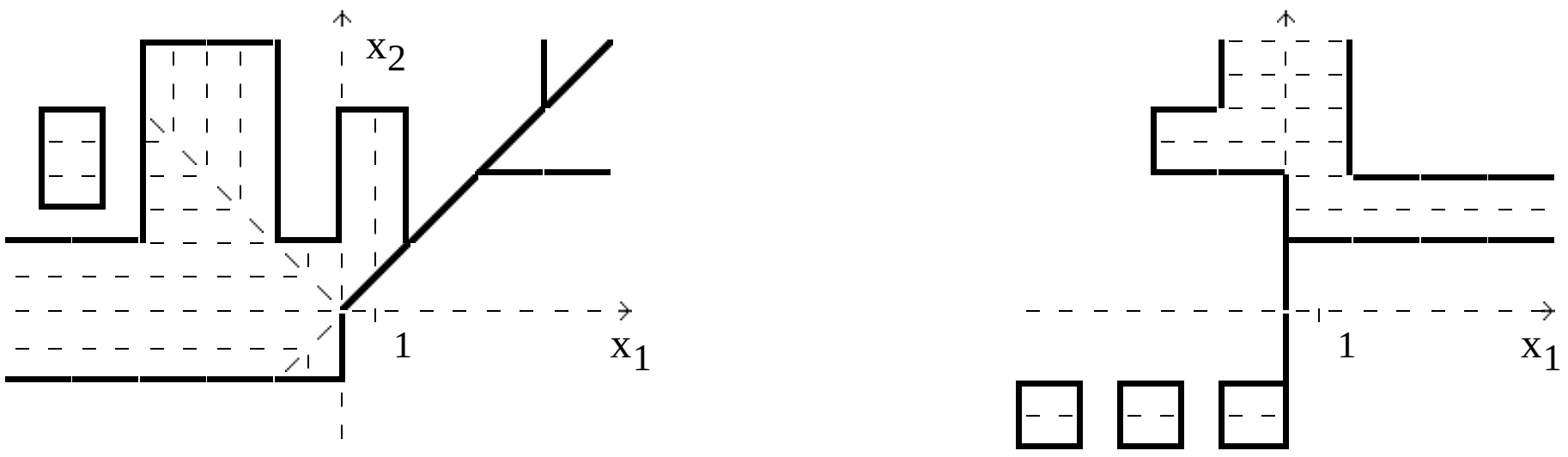

On the other hand, a subspace $X \subset E^2$ (as at the right hand, above) is *telescopic with respect to the axis* $x_2$ if the homotopy $\gamma = \alpha{\times}0_{id}$ (3.4.5) satisfies a similar condition: for any $x = (x_1, x_2) \in X$, the path $\gamma(-, x)$, from $(0, x_2)$ to $x$, is contained in $X$. Then, $X$ admits as a positive deformation retract its intersection with the axis $x_2$, and is contractible iff this intersection is so.

One defines similarly telescopic domains in $\mathbf{R}^n$ or $\mathbf{Z}^n$, with respect to any point, or more generally to any intersection of hyperplanes $x_i$ = constant.

**3.6. Variable jumps.** These homotopies are easily adapted to the $t_\varepsilon$-structure of $E^n$. Less trivially, to prove that the c-space $t_\varepsilon X \subset t_\varepsilon \mathbf{R}^2$ described in figure 1.8.1 is contractible for $\varepsilon \geq 3$, we need a generalised telescopic homotopy centred at the horizontal axis, with "variable vertical jumps" 1, 3, 1 (all $\leq \varepsilon$, and adjusted to jump over the singularity $A \cup B \cup C$).

To make this precise, let us consider the real or integral line $t_\varepsilon E$, with the $t_\varepsilon$-structure ($\varepsilon > 0$). We have now a positive (generalised) *telescopic homotopy* with variable jumps in $t_\varepsilon E$, centred at $a_0$

(1) $\alpha$: $p \to id$: $t_\varepsilon E \to t_\varepsilon E$,

$\alpha(i, x) = a_0 \vee (a_i \wedge x) \quad (x \geq a_0)$, $\qquad \alpha(i, x) = a_0 \wedge (a_{-i} \vee x) \quad (x \leq a_0)$,

between the constant map $p(x) = a_0$ and the identity, determined by its *characteristic sequence* $(a_i)$, an *unbounded* increasing combinatorial mapping $a: \mathbf{Z} \to t_\varepsilon E$ (no lower nor upper bound); also determined by its *centre* $a_0 \in E$ and the *sequence of jumps* $s_i = a_i - a_{i-1} \in E$ $(0 \leq s_i \leq \varepsilon)$. The homotopy (1) reduces to the standard telescopic homotopy (centred at the origin) when $\varepsilon = 1$ and $a_i = i$.

To verify that $\alpha: \mathbf{Z} \times t_\varepsilon E \to t_\varepsilon E$ is indeed a map, after invoking the lattice structure of $t_\varepsilon E$ (1.7.1), it is sufficient to consider two linked pairs on opposite sides of $a_0$

(2) $(i, x),\ (i', x') \in \mathbf{Z} \times t_\varepsilon E, \qquad x' < a_0 < x, \quad x - x' \leq \varepsilon, \quad |i - i'| \leq 1,$

and note that $0 \leq \alpha(i, x) - \alpha(i', x') \leq x - x' \leq \varepsilon$.

Extending 3.4, we have (generalised) telescopic homotopies in $t_\varepsilon \mathbf{R}^n$ or $t_\varepsilon \mathbf{Z}^n$, with respect to any intersection of hyperplanes $x_i$ = constant, determined by unbounded sequences $\mathbf{Z} \to t_\varepsilon E$ (one for each variable on which our homotopy is not trivial). Extending 3.5, we have (generalised) telescopic domains in $t_\varepsilon \mathbf{R}^n$ or $t_\varepsilon \mathbf{Z}^n$. Coming back to the initial example, and letting $\varepsilon \geq 3$, the simplicial complex $t_\varepsilon X \subset t_\varepsilon \mathbf{R}^2$ is telescopic with respect to the horizontal axis, with jumps $s_1 = 1$, $s_2 = 3$, $s_3 = 1$ on the second variable (arbitrarily completed), and homotopy equivalent to its intersection with the vertical axis, which is contractible.

**3.7. Grid retracts.** $\mathbf{Z}$ is a positive deformation retract of $\mathbf{R}$, with retraction the integral-part map $[-]: \mathbf{R} \to \mathbf{Z}$. In fact, the discontinuous map $[-]$ is combinatorial because $[x] - [y] < x - y + 1$, so that $[x] - [y] \geq 2$ implies $x - y > 1$. We have thus a positive homotopy (for the inclusion $j: \mathbf{Z} \subset \mathbf{R}$)

(1) $\alpha: j[-] \to \mathrm{id}: \mathbf{R} \to \mathbf{R}, \qquad \alpha(i, x) = [x] \vee (i \wedge x).$

Note that each path $\alpha(-, x)$ has an image reduced to two points, $\{[x], x\}$ and admits a support reduced to two indices, $i = [x], [x]+1$; but $\mathbf{Z}$ is not an immediate deformation retract of $\mathbf{R}$, nor even a bounded one, as showed in 3.4.

Therefore, $\mathbf{Z}^n$ is a positive retract of $\mathbf{R}^n$. More generally, let us fix $\varepsilon > 0$. Then, the *n-dimensional* ε*-grid* $t_\varepsilon(\varepsilon\mathbf{Z}^n)$ is a positive deformation retract of $t_\varepsilon(\mathbf{R}^n)$, with retraction the n-th power of $[x]_\varepsilon = \varepsilon.[\varepsilon^{-1}.x]$,

(2) $[-]_\varepsilon: \mathbf{R} \to \varepsilon\mathbf{Z}, \qquad x - \varepsilon < [x]_\varepsilon \leq x.$

Now, if a metric subspace $X \subset \mathbf{R}^n$ contains its ε*-grid* $[X]_\varepsilon$, then $t_\varepsilon[X]_\varepsilon$ is a positive deformation retract of $t_\varepsilon X$, the ε*-grid retract*. Thus, the two (connected!) combinatorial subspaces $t_\varepsilon X$, $t_\varepsilon Y \subset t_\varepsilon \mathbf{R}^2$ displayed below are homotopy equivalent: $t_\varepsilon X \simeq t_\varepsilon[X]_\varepsilon = t_\varepsilon[Y]_\varepsilon \simeq t_\varepsilon Y$

(3)

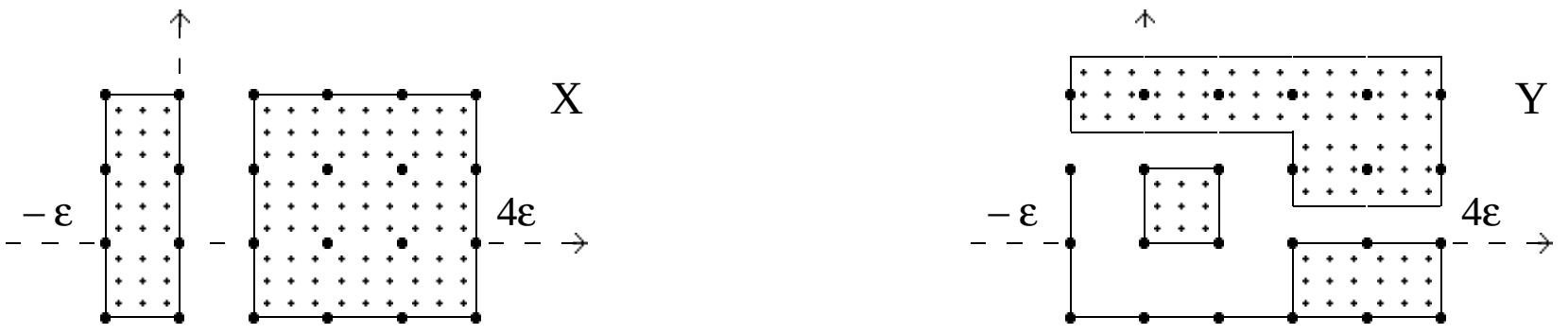

## 4. Homotopy pullbacks and the fibre sequence

The problems which arise in studying combinatorial homotopy will result more clearly by comparison with a general frame holding in classical situations, like topological spaces and chain complexes. We review some general notions on "categories with homotopies", within a setting developed in [Gr2]. **Cs** and **Cs***  have a defective structure of this type, which will be extended in the next Section.

**4.1. Abstract homotopy structures.** An h4-*category* **A** is a sort of relaxed 2-category, abstracting the "nearly 2-categorical" properties of **Top** and other situations where homotopy arose.

Precisely, the category **A** is equipped with *cells* or *homotopies* $\alpha\colon f \to g\colon X \to Y$ and the following "operations" of various arities, on maps and cells (the "vertical structure" is written in additive notation; the whisker composition will also be written by juxtaposition)

(a) *whisker composition* (or reduced horizontal composition) of maps and cells:

$v{\circ}\alpha{\circ}u\colon\ vfu \to vgu\colon\ X' \to Y'$ (for $u\colon X' \to X,\ v\colon Y \to Y'$),

(b) *trivial homotopies* (or vertical identities): $0_f\colon f \to f$,

(c) *reversion* (or vertical involution): $-\alpha\colon g \to f$ (for $\alpha\colon f \to g$),

(d) *concatenation* (or vertical composition): $\alpha + \beta\colon f \to h$ (for $\alpha\colon f \to g,\ \beta\colon g \to h$),

(e) the *2-homotopy relation* $\alpha \simeq_2 \alpha'$, an equivalence relation for parallel homotopies $\alpha, \alpha'\colon f \to g$,

satisfying the following axioms (the last property, (hc4c), is the *weak reduced interchange*):

(hc0) $1_B{\circ}\alpha{\circ}1_A = \alpha$, $v{\circ}0_f{\circ}u = 0_{vfu}$, $v'{\circ}(v{\circ}\alpha{\circ}u){\circ}u' = (v'v){\circ}\alpha{\circ}(uu')$,

(hc1) $-0_f = 0_f$, $-(-\alpha) = \alpha$, $-(v{\circ}\alpha{\circ}u) = v{\circ}(-\alpha){\circ}u$,

(hc2) $0_f + 0_f = 0_f$, $v{\circ}(\alpha+\beta){\circ}u = (v{\circ}\alpha{\circ}u) + (v{\circ}\beta{\circ}u)$,

(hc3) $-(\alpha+\beta) = (-\beta) + (-\alpha)$,

(hc4a) $v{\circ}\alpha{\circ}u \simeq_2 v{\circ}\alpha'{\circ}u$, $-\alpha \simeq_2 -\alpha'$,
$\alpha+\beta \simeq_2 \alpha'+\beta'$ (for $\alpha \simeq_2 \alpha'$, $\beta \simeq_2 \beta'$),

(hc4b) $0_f+\alpha \simeq_2 \alpha \simeq_2 \alpha+0_g$, $\alpha+(-\alpha) \simeq_2 0_f$,
$-\alpha+\alpha \simeq_2 0_g$, $(\alpha+\beta)+\gamma \simeq_2 \alpha+(\beta+\gamma)$,

(hc4c) $u{\circ}\alpha + \varphi{\circ}g \simeq_2 \varphi{\circ}f + v{\circ}\alpha$ (for $\alpha\colon f \to g\colon X \to Y$, $\varphi\colon u \to v\colon Y \to Z$).

Weaker notions are also used. Thus, an *h-category* (or h0-category) has just whisker composition and trivial homotopies, satisfying (hc0); an *h1-category* has also a reversion, under (hc1); an *h2-category* is an h-category with concatenation, satisfying (hc2); finally, an *h3-category* has both reversion and concatenation, under (hc0-3). An *hi-functor* $F\colon \mathbf{A} \to \mathbf{B}$ between hi-categories (i = 0,... 4) preserves their whole structure; an *hi-transformation* $u\colon F \to G\colon \mathbf{A} \to \mathbf{B}$ of hi-functors is required to be 2-natural ($uB{\circ}F\alpha = G\alpha{\circ}uA$).

In an h1-category, the homotopy relation $f \simeq g$ (there exists a homotopy $f \to g$) is reflexive, symmetric and *weakly* consistent with composition: $f \simeq g$ implies $vfu \simeq vgu$. The map $f\colon X \to Y$ is a *homotopy equivalence* if it has a homotopy inverse g ($gf \simeq 1$, $fg \simeq 1$). Two objects are *homotopy equivalent* if there is a finite sequence of homotopy equivalences connecting them. The map f is a *fibration* if it satisfies the usual lifting property of homotopies: for every $x: A \to X$, every homotopy $\psi: fx \to y: A \to Y$ lifts to some homotopy $\varphi: x \to x'$, with $f\varphi = \psi$ (and $fx' = y$). If **A** is h3, all this can be simplified: the homotopy relation $f \simeq g$ is a congruence; its quotient Ho**A** = **A**/$\simeq$ is called the *homotopy category* of **A**; f is a homotopy equivalence iff its homotopy class [f] is an isomorphism of Ho**A**.

A functor F: **A** → **B** between h1-categories is *homotopy invariant* if it preserves the homotopy relation; then it also preserves homotopy equivalences and fibrations. This also makes sense for an ordinary category **B**, viewed as a trivial h1-category (with formal trivial homotopies); then, a functor F: **A** → **B** is homotopy invariant iff it identifies all pairs of homotopic maps.

A *strict* h4-category is an h4-category whose 2-homotopy relation $\simeq_2$ is the equality; it is not difficult to see that this is equivalent to a 2-category whose cells are invertible, i.e. a *groupoid-enriched category* ([Gr2], 1.4). Every h4-category has an associated strict one, $\mathrm{Ho}_2$**A** = **A**/$\simeq_2$, consisting of the same objects, same maps and *tracks* (homotopies modulo $\simeq_2$).

More generally, an h4-category will be said to be *semiregular* if its concatenation is strictly associative, with strict identities; it is then a sesqui-category in the sense of Street [St] (every set of maps and homotopies **A**(X, Y) is a category, coherently with left and right composition *with maps*), with additional structure: reversion and the 2-homotopy relation. **A** is *regular* if, moreover, the reversion gives strict inverses.

The 2-category $\mathbf{Cat}_i$ of small categories, functors and natural isomorphisms is thus *strict h4*, as well as **Gpd** (small groupoids). The category **Dm** of chain complexes, with the usual homotopies, is *regular h4*, while **Top** (topological spaces) and **Top*** (pointed spaces) are just *h4-categories* [Gr2]. The category **Da** of chain algebras, with multiplicative homotopies, is an h-category having a forgetful functor into **Dm**, an h4-category; we need to work with both structures [Gr3]. Similarly **Cs** and **Cs*** *are only h1-categories* (3.1), and will be embedded in richer structures (Section 5).

In an h-category, a 2-*terminal* object $\top$ has to satisfy a 2-dimensional universal property, implying the usual 1-dimensional one: for every object X there is precisely one homotopy $X \Rightarrow \top$, i.e. one map $X \to \top$ and one endocell of the latter, its vertical identity. Similarly, one defines a 2-*product* $X = \Pi X_i$ (for each family of homotopies $\alpha_i: A \Rightarrow X_i$ there is precisely one homotopy $\alpha: A \Rightarrow X$ whose projections are the given ones: $p_i \circ \alpha = \alpha_i$). In **Cs** and **Cs***, *finite* products are 2-products, because P preserves them (2.4), and arbitrary sums are 2-sums (as it trivially happens whenever homotopies are represented by a path functor). All the other h-categories considered above have arbitrary 2-products and 2-sums (as they possess both a path and a cylinder functor).

**4.2. Homotopy pullbacks.** Let **A** be an h-category. The *standard homotopy pullback*, or *h-pullback*, of two maps f, g having the same codomain, is defined as an object P(f, g) equipped with two maps p, q and a homotopy $\xi: fp \to gq: P(f, g) \to Y$, as in the left-hand diagram below

$$\text{(1)}\qquad \begin{array}{ccc} X & \xrightarrow{f} & Y \\ {\scriptstyle p}\uparrow & \searrow{\scriptstyle \xi} & \uparrow{\scriptstyle g} \\ P(f, g) & \xrightarrow[q]{} & Z \end{array} \qquad\qquad \begin{array}{ccc} Y & \xrightarrow{1} & Y \\ {\scriptstyle \partial^-}\uparrow & \searrow{\scriptstyle \delta} & \uparrow{\scriptstyle 1} \\ PY & \xrightarrow[\partial^+]{} & Y \end{array}$$

satisfying the obvious universal property (of comma squares): given similar data $\xi'\colon fp' \to gq'\colon W \to Y$ there is precisely one map $u\colon W \to P(f, g)$ such that $pu = p'$, $qu = q'$, $\xi u = \xi'$. The solution $(P(f, g); p, q, \xi)$ is determined up to isomorphism.

The triple $(p, q; \xi)$ is jointly monic. The h-pullback $PY = P(1_Y, 1_Y)$, as in the right hand diagram above, is called the *path-object* of $Y$. It comes equipped with a *universal* homotopy $\delta\colon \partial^- \to \partial^+\colon PY \to Y$ (*path evaluation*), establishing a bijection between maps $\alpha\colon W \to PY$ and homotopies $\delta\circ\alpha\colon \partial^-\alpha \to \partial^+\alpha\colon W \to Y$. Plainly, *an h-category with path-objects* amounts to *a category equipped with an endofunctor* $P$ *and three natural transformations* $\partial^-, \partial^+\colon P \to 1$ (*faces*), $e\colon 1 \to P$ (*degeneracy*), satisfying $\partial^- e = 1 = \partial^+ e$ (dually to Kan's setting for a cylinder functor [Ka]). If **A** is h1, the (existing) h-pullbacks are symmetric: $(P(f, g); q, p, -\xi)$ is the h-pullback of $g, f$.

Homotopy pullbacks can be calculated by means of the path object $PY$ and of ordinary pullbacks, when such constructs exist. In fact, the general h-pullback $P(f, g)$ can be obtained as the ordinary limit of the diagram (2), which amounts to two ordinary pullbacks

$$\text{(2)}\qquad X \xrightarrow{f} Y \xleftarrow{\partial^-} PY \xrightarrow{\partial^+} Y \xleftarrow{g} Z$$

In **Top**, homotopies are represented by the path space $PY = Y^{[0, 1]}$ (with the compact-open topology), and by its left adjoint, the cylinder functor $IX = [0, 1]\times X$. Accordingly, the h-pullback of $f, g$ as in (1) is $P(f, g) = \{(x, \beta, z) \in X\times PY\times Z \mid \beta(0) = f(x),\ \beta(1) = g(z)\}$, while the h-pushout $I(f, g)$ of $f\colon X \to Y$, $g\colon X \to Z$ is realised pasting the spaces $Y$ and $Z$ over the bases of the cylinder $IX = [0, 1]\times X$, along $f$ and $g$. Homotopy pullbacks and pushouts in **Dm** and **Da** are also well known (cf. [Gr3]). In **Cs** and **Cs***, the h-pullback is similarly obtained as a limit (4.5), whereas h-pushouts do not exist, in general: $P$ has no left adjoint (2.4).

**4.3. The fibre sequence.** Let **A** be a *pointed* h-category with h-pullbacks, as **Top***, **Dm**, **Da** and **Cs***: it has a zero object $0$, which is 2-initial and 2-terminal (4.1). The (upper) *homotopy kernel*, or h-kernel, or homotopy fibre $Kf = K^+f = P(0, f)$ is the h-pullback of $0 \to Y$, along $f$

$$\text{(1)}\qquad \begin{array}{ccc} X & \xrightarrow{f} & Y \\ {\scriptstyle x}\uparrow & \nwarrow{\scriptstyle \xi} & \uparrow \\ Kf & \longrightarrow & 0 \end{array} \qquad\qquad \mathrm{hkr}(f) \ = \ (0, x; \xi)$$

the map $x\colon Kf \to X$ will be called the *h-kernel map* of $f$, and written $\mathbf{k}(f)$; the homotopy $\xi\colon 0 \to fx\colon Kf \to Y$ will be called the *h-kernel homotopy* of $f$. In particular, homotopy kernels also yield the (upper) *cocone* endofunctor $K$ and the *loop* endofunctor $\Omega$, already described in 2.4.3

(2) $K\colon \mathbf{A} \to \mathbf{A}$, $\qquad KX \ = \ K^+X \ = \ K(1_X)$,

(3) $\Omega\colon \mathbf{A} \to \mathbf{A}$, $\qquad K(0 \to X)$,

the latter equipped with a universal homotopy $\omega\colon 0 \to 0\colon \Omega X \to X$ (*loop evaluation of* X).

Extending the standard situation of pointed topological spaces, every map $f\colon X \to Y$ has a natural *fibre sequence*, or *dual Puppe sequence* (cf. [Pu])

$$(4) \qquad \ldots \longrightarrow \Omega Kf \xrightarrow{\Omega x} \Omega X \xrightarrow{\Omega f} \Omega Y \xrightarrow{d} Kf \xrightarrow{x} X \xrightarrow{f} Y$$

where $x = \mathbf{k}(f)$ is the h-kernel map of f, while the *differential* $d\colon \Omega Y \to Kf$ is determined by the following conditions (for $hkr(f) = (0, x; \xi)$)

$$(5) \quad xd = 0, \qquad\qquad \xi \circ d = \omega_Y.$$

The fibre sequence of f can be linked to the one of its h-kernel map x by the *comparison map* $u\colon \Omega Y \to Kx$, forming a diagram which is commutative, *except for the left-hand* comparison *square*

$$(6) \qquad \begin{array}{ccccccccc} \ldots\ \Omega X & \xrightarrow{\Omega f} & \Omega Y & \xrightarrow{d} & Kf & \xrightarrow{x} & X & \xrightarrow{f} & Y \\ \| & \# & \downarrow u & & \| & & \| & & \\ \ldots\ \Omega X & \xrightarrow[d']{} & Kx & \xrightarrow[y]{} & Kf & \xrightarrow[x]{} & X & & \end{array}$$

$$(7) \quad yu = d, \qquad \eta \circ u = 0\colon xd \to 0\colon \Omega Y \to X \qquad (hkr(x) = (0, y; \eta)).$$

Taking on this procedure, we end by linking the fibre sequence of f to the sequence of all its iterated h-kernels ($x_n = \mathbf{k}(x_{n-1})$, with $x_0 = f$), in a *fibre diagram* ([Gr2], 5.5; 6.4), generally non-commutative, whose columns consist of comparison maps and their images by powers $\Omega^n$.

All this behaves well *when* **A** *is pointed h4 and its h-kernels satisfy also a 2-dimensional universal property* (4.4). Then, as proved in [Gr2] (5.6; 6.6), the following results hold.

*Regularity properties.* All h-kernel maps are fibrations; all comparison maps (as well as their $\Omega^n$-images) are homotopy equivalences; all comparison squares are homotopy anti-commutative (with respect to the reversion in loop-objects). A loop-object $\Omega X$ has a natural structure of internal group "up to homotopy", by loop-concatenation; the two induced structures over $\Omega^2 X$ have the same identity and satisfy the middle-interchange property, up to homotopy.

It follows that each object S defines a sequence of homotopy-invariant functors $\pi_n = \pi_n^S$

$$(8) \quad \pi_n = [S, \Omega^n(-)]\colon \mathbf{A} \to \mathbf{Set}^*, \qquad (\pi_n = \pi_0 \Omega^n),$$

taking values in groups for $n \geq 1$, and abelian groups for $n \geq 2$. Any map f has an exact sequence of pointed mappings

$$(9) \qquad \ldots \longrightarrow \pi_1 Kf \xrightarrow{\pi_1 x} \pi_1 X \xrightarrow{\pi_1 f} \pi_1 Y \xrightarrow{\pi_0 d} \pi_0 Kf \xrightarrow{\pi_0 x} \pi_0 X \xrightarrow{\pi_0 f} \pi_0 Y.$$

**4.4. The 2-dimensional property.** If **A** is an h4-category, the h-pullback $(P(f, g); p, q; \xi)$ is said to be *regular* (or an h4-pullback) if, given two maps $a, b\colon W \to P(f, g)$, as in the left hand diagram below, and two homotopies $\varphi\colon pa \to pb\colon W \to X$, $\psi\colon qa \to qb\colon W \to Y$ coherent with $\xi$ (i.e., $f\varphi + \xi b \simeq_2 \xi a + g\psi$, as in the right-hand diagram of homotopies)

$$(1)\qquad W \underset{b}{\overset{a}{\rightrightarrows}} P(f, g) \qquad \begin{array}{ccccc} & & X & & \\ & {}^{p}\nearrow & & \searrow^{f} & \\ P(f,g) & & \downarrow \xi & & Z \\ & {}_{q}\searrow & & \nearrow_{g} & \\ & & Y & & \end{array} \qquad\qquad \begin{array}{ccc} gqa & \xrightarrow{g\psi} & gqb \\ \xi a \uparrow & \simeq_2 & \uparrow \xi b \\ fpa & \xrightarrow[f\varphi]{} & fpb \end{array}$$

there is *some* homotopy $\alpha: a \to b$ which lifts $\varphi$ and $\psi$: $p\alpha = \varphi$, $q\alpha = \psi$.

The basic properties of *regular* h-pullbacks, as homotopy invariance, pasting, and reflection of homotopy equivalences, are studied in [Gr2] (Sections 2, 3), together with their consequences on the fibre sequence of a map, in the pointed case (Section 6). **Top**, **Top*** and **Dm** have regular h-pullbacks (and h-pushouts).

**Da** does not have a concatenation of homotopies, and we cannot even formulate the coherence hypothesis in (1); but its h-pullbacks are preserved by the forgetful functor to **Dm**, where they satisfy the 2-dimensional property; we can thus draw the fibre sequence of a map within **Da**, and study its 2-dimensional properties in **Dm** [Gr3]. Similarly, **Cs*** is just an h1-category; but its h-pullbacks are preserved by the embedding into Ps**Cs*** where they are regular (Section 5); we shall use this to study the fibre sequence of a map and the homotopy groups (Section 6).

**4.5. Homotopy pullbacks of simplicial complexes.** In **Cs**, the h-pullback P(f, g) is obtained as in 4.2.2, by means of the path object PY and of two ordinary pullbacks

(1) $P(f, g) \;=\; \{(x, a, z) \in X{\times}PY{\times}Z \mid \partial^-a = f(x),\ \partial^+a = g(z)\} \;\subset\; X{\times}PY{\times}Z.$

We deduce now that P preserves h-pullbacks (via its extension to homotopies, 3.1.2). In fact, P preserves finite limits (2.4), in particular the ordinary limit 4.2.2; modifying it via the automorphism $s: P^2Y \to P^2Y$

$$(2)\qquad \begin{array}{ccccccccc} PX & \xrightarrow{Pf} & PY & \xleftarrow{P\partial^-} & P^2Y & \xrightarrow{P\partial^+} & PY & \xleftarrow{Pg} & PZ \\ \| & & \| & & s\uparrow & & \| & & \| \\ PX & \xrightarrow{Pf} & PY & \xleftarrow{\partial^-P} & P^2Y & \xrightarrow{\partial^+P} & PY & \xleftarrow{Pg} & PZ \end{array}$$

it is clear that $(P(P(f, g)), Pp, Pq, \hat{P}\xi)$ is the h-pullback of Pf, Pg.

This fact produces a 2-dimensional property of P(f, g) for *double* homotopies, which is not sufficient to prove the "regularity properties" considered above (4.3), but will be of use later (5.7). Given two homotopies $\varphi: pa^- \to pa^+: W \to X$, $\psi: qa^- \to qa^+: W \to Z$

$$(3)\qquad W \underset{a^+}{\overset{a^-}{\rightrightarrows}} P(f, g) \qquad \begin{array}{ccccc} & & X & & \\ & {}^{p}\nearrow & & \searrow^{f} & \\ P(f,g) & & \downarrow \xi & & Y \\ & {}_{q}\searrow & & \nearrow_{g} & \\ & & Z & & \end{array} \qquad\qquad \begin{array}{ccc} gqa^- & \xrightarrow{Pg.\psi} & gqa^+ \\ \xi a^- \uparrow & \Phi & \uparrow \xi a^+ \\ fpa^- & \xrightarrow[Pf.\varphi]{} & fpa^+ \end{array}$$

coherent with $\xi$ *by means of a double homotopy* $\Phi: W \to P^2Y$ with faces as above, there is *some* homotopy $\alpha: a^- \to a^+$ which lifts $\varphi$ and $\psi$: $p{\circ}\alpha = \varphi$, $q{\circ}\alpha = \psi$.

The proof is an easy application of the universal property of the h-pullback $P(P(f, g))$: we have two maps $\varphi: W \to PX$, $\psi: W \to PZ$ and a homotopy $\Phi: W \to P(PY)$ such that

(4) $\partial^{-}P.\Phi = Pf.\varphi, \qquad \partial^{+}P.\Phi = Pg.\psi,$

whence a unique map $\alpha: W \to P(P(f, g))$ such that

(5) $p\circ\alpha = Pp.\alpha = \varphi, \qquad q\circ\alpha = Pq.\alpha = \psi, \qquad (\hat{P}\xi)\circ\alpha = \Phi.$

Finally, the faces of $\alpha$ are indeed $a^{\kappa}$, as detected by the jointly monic triple $(p, q, \xi)$:

(6) $p(\partial^{\kappa}\alpha) = \partial^{\kappa}.Pp.\alpha = \partial^{\kappa}\varphi = pa^{\kappa}, \qquad q(\partial^{\kappa}\alpha) = qa^{\kappa},$

$\xi(\partial^{\kappa}\alpha) = \partial^{\kappa}P.P\xi.\alpha = \partial^{\kappa}P.s.\hat{P}\xi.\alpha = P\partial^{\kappa}.\Phi = \xi.a^{\kappa}.$

**4.6. Other path functors.** We end this section by examining three other possible "path functors" for **Cs**, their "homotopies" and their disadvantages with respect to our previous functor P. Such functors, including P, are linked by natural transformations, componentwise surjective

(1) $\mathtt{P}X \to \mathtt{P}'X \to PX \to \mathbf{P}X.$

(a) The c-space $\mathtt{P}X \subset DJ\times PX$ of *paths with duration* consists of all pairs $\hat{a} = (\rho, a)$ where $\rho\in J$ is an admissible support of $a\in PX$ and DJ is discrete. (Equivalently, $\mathtt{P}X = \Sigma_{\rho}\, \{\rho\}\times X^{[\rho^{-}, \rho^{+}]}$.) This allows one to concatenate paths in a natural way, preserved by maps. Formally, we get "nearly" a cubical comonad with symmetries and concatenation (the only axiom which fails is absorbency, cf. 2.4.2). A $\mathtt{P}$-*homotopy* $\hat{\alpha} = (\rho, \alpha): X \to \mathtt{P}Y$ amounts to an ordinary homotopy $\alpha$ and a support-map $\rho: X \to DJ$ for it (a mapping constant on connected components). Thus, two maps are $\mathtt{P}$-homotopic iff they are linked by an ordinary homotopy bounded on connected components, $f \simeq_{c} g$ (3.2); but a $\mathtt{P}$-homotopy also carries a particular $\rho$. Various defects arise: first, the terminal object $\{*\}$ is not 2-terminal (the endo-$\mathtt{P}$-homotopies of $X \to \{*\}$ correspond to $\rho\in J$); similarly, h-pullbacks do not satisfy the 2-dimensional property, but only "do so up to modifying supports", not only in $\mathtt{P}$-homotopies but also in *maps* with values in the h-pullback $\mathtt{P}(f, g)$; finally, and more concretely, the standard line **Z** is not contractible by $\mathtt{P}$-homotopies (as for bounded homotopies, 3.2).

(b) Taking $\mathtt{P}'X \subset \mathbf{J}\times PX$ (with the contiguity structure on $J \subset \mathbf{Z}^2$), the last defect is amended but we miss concatenation of homotopies (and interchange), loosing any advantage with respect to P.

(c) Finally, the quotient modulo delays $\mathbf{P}X = (PX)/\!\equiv$ (2.6) looks promising at first, since **P**-homotopies form a semiregular h4-category. But there seems to be no 2-dimensional property for h-pullbacks $\mathbf{P}(f, g)$, essentially because **P** does not preserve binary products (nor pullbacks, as it preserves the terminal). In fact, if $Y = X_1\times X_2$, the congruence $\equiv$ of PY implies the product of the congruences of $PX_k$ (if $a_k.d = b_k.e$ for $k = 1, 2$, then $a_k \equiv b_k$); but the converse is false: take two different paths $a \equiv a': \mathbf{Z} \to X$, and consider the paths $\langle a, a\rangle, \langle a, a'\rangle: \mathbf{Z} \to X\times X$; they have congruent projections on X, but are not congruent: $\langle a, a\rangle.d = \langle a, a'\rangle.e$ would give $ad = ae$, $ad = a'e$, whence $ae = a'e$ and $a = a'$, as all delays are surjective.

## 5. Bounded homotopies and pseudo-homotopies

Bounded homotopies have a concatenation, and satisfy the axioms of h4-categories in a slightly laxified form (5.1). But they do not include the structural homotopies of h-pullbacks, and we need to embed **Cs** in a larger lax h4-category, Ps**Cs**.

**5.1. Lax h4-categories.** The structures we are going to construct satisfy the axioms of h4-category (4.1) in a laxified form, already used in [Gr5] for the homotopy structure of Stasheff's $A_\infty$-algebras (strongly homotopy associative differential algebras). All the main results of [Gr2] are still valid for this extension, with minimal modifications in proofs. (Moreover, this laxified version is probably more natural than the original one, and might replace it.)

A *lax h4-category* **A** is an h-category equipped with the same structure considered above (in 4.1: reversion $-\alpha$; concatenation $\alpha + \beta$; the equivalence relation $\alpha \simeq_2 \alpha'$), under the only axiom (hc4). The latter can be summarised saying that all operations ($v \circ \alpha \circ u$; $-\alpha$; $\alpha + \beta$) are consistent with $\simeq_2$ and make the quotient $\mathrm{Ho}_2\mathbf{A} = \mathbf{A}/\simeq_2$ into a 2-category with invertible cells.

As a consequence, the other axioms (hc.1-3) hold strictly in the quotient; but in **A** they only have to hold up to 2-homotopy. In fact, the only modification which really takes place *here* concerns (hc2), which will only hold in laxified form, $v \circ (\alpha+\beta) \circ u \simeq_2 (v \circ \alpha \circ u) + (v \circ \beta \circ u)$, essentially because of the *non-natural* concatenation we use here (as for paths; cf. 2.5.3).

**5.2. The concatenation of bounded homotopies.** We begin by studying the structure of $\mathbf{Cs}_b$, the category of simplicial complexes equipped with *bounded* homotopies (admitting a constant support). First, the operations considered in 3.1 (whisker composition, trivial homotopies, reversion) restrict to bounded homotopies, giving raise to an h1-structure.

Now, the concatenation of paths (2.5) can be transferred to bounded homotopies, which are in a sense an equivalent notion (by the cartesian closedness of **Cs**, 1.4). Trivially, a path is a (bounded) homotopy $a: x \to x': \{*\} \to X$ of its end points. On the other hand, a *bounded* homotopy $\alpha: \mathbf{Z}\times X \to Y$ can be viewed as *a path in* $\mathrm{Hom}(X, Y) = Y^X$ (a line $\alpha: \mathbf{Z} \to Y^X$ having a support).

Thus, a bounded homotopy $\alpha: X \to PY$ has already been given a *standard support* $\rho(\alpha) \in J$; two consecutive bounded homotopies $\alpha: f \to g$, $\beta: g \to h$ have a (standard) *concatenation* $\gamma = \alpha + \beta$

(1) $\gamma(i, x) = \alpha(i - \rho^-(\beta), x)$, for $i \leq \rho^+(\alpha) + \rho^-(\beta)$,

$\gamma(i, x) = \beta(i - \rho^+(\alpha), x)$, for $i \geq \rho^+(\alpha) + \rho^-(\beta)$,

and various *admissible concatenations*. Two bounded homotopies $\alpha, \alpha': f \to g$ are *congruent* (up to delays), $\alpha \equiv \alpha'$, when there are delays $d, d': \mathbf{Z} \to \mathbf{Z}$ (2.6) such that $\alpha.(d\times X) = \alpha'.(d'\times X): \mathbf{Z}\times X \to Y$. All admissible concatenations of two given bounded homotopies are congruent.

Occasionally, and technically, we use more general admissible concatenations, for a *right bounded* $\alpha$ and a *left bounded* $\beta$ (3.1). In fact, the procedure (1) makes sense whenever the support of $\alpha$ has a (constant) *right bound* $\rho^+(\alpha)$, and the one of $\beta$ a left bound $\sigma^-(\beta)$; moreover, if $\alpha$ (resp. $\beta$) is bounded, the result is left (resp. right) bounded. For instance, this allows us to prove a fact anticipated in 3.3: if the c-space $X$ has a *finite positive* deformation retract $u: S \to X$, then two maps $f, g: X \to Y$ are homotopic iff their restrictions $fu, hu: S \to Y$ are so. Indeed, if $fu \simeq gu: S \to Y$, we have a *bounded* homotopy $\varphi: fup \to gup: X \to Y$, which we can concatenate with the *negative homotopy* $-f\alpha$: $f \to fup$ (at the left) and the *positive homotopy* $g\alpha$: $gup \to g$ (at the right), forming a homotopy $f \to g$. (One can also extend concatenation to homotopies bounded on connected components, using supports $\rho$: $X \to DJ$; this is not used here.)

**5.3. Bounded double homotopies.** Similarly, a *bounded* double homotopy $\Phi$: $X \to P^2Y$ (3.1) is a double path in $Y^X$, and a *bounded* 2-homotopy is a 2-path in $Y^X$; all previous constructions on 2-dimensional paths (2.7-8) can be transferred. Our selection determines a standard support $(\rho(\Phi), \sigma(\Phi)) \in J^2$ consistent with interchange; it also yields the (standard) vertical pasting of vertically consecutive bounded double homotopies, and – symmetrically – the horizontal analogue

(1) $s(\Phi +_v \Psi) = s\Phi +_h s\Psi$ $\qquad (\partial^+P.\Phi = \partial^-P.\Psi)$,

with the following non-trivial relations on the concatenation of faces

(2) $P\partial^\kappa.(\Phi +_v \Psi) \equiv P\partial^\kappa.\Phi + P\partial^\kappa.\Psi$, $\qquad \partial^\kappa P.(\Phi +_h \Psi) \equiv \partial^\kappa P.\Phi + \partial^\kappa P.\Psi$.

Both operations are associative, satisfy the cancellation property and, together, the four middle interchange law up to 2-dimensional delays (2.7.6). Both have identities and involution

(3) $0_v(\alpha) = eP.\alpha$, $\qquad 0_h(\alpha) = Pe.\alpha$ $\qquad$ (*vertical and horizontal identity of* $\alpha$),

(4) $-_v \Phi = rP.\Phi$ $\qquad -_h \Phi = Pr.\Phi$ $\qquad$ (*vertical and horizontal reversion of* $\Phi$).

A *bounded 2-homotopy* $\Phi$: $\alpha \simeq_2 \beta$: $f \to g$ is a bounded double homotopy whose horizontal faces are $\alpha, \beta$, while its vertical faces are trivial $(0_f, 0_g)$. Such double homotopies are stable under horizontal and vertical reversion or pasting. The *bounded 2-homotopy relation* $\alpha \simeq_{2b} \beta$ (also written $\alpha \simeq_2 \beta$, within $\mathbf{Cs}_b$) is an equivalence relation. The relation $\alpha \equiv \beta$ implies $\alpha \simeq_2 \beta$ (*Caterpillar Lemma*), so that all admissible concatenations of a pair of bounded homotopies are 2-homotopic.

**5.4. Theorem: The structure of bounded homotopies.** $\mathbf{Cs}_b$ is a *lax* h4-category (5.1) with semiregular concatenation. Moreover, the axioms (hc1, 3) are satisfied.

**Proof.** All homotopies and double homotopies used here are bounded. The proof is as in [Gr4], thm. 2.9, with some modifications due to semiregularity (which simplifies things) and the lax relation in 5.3.2 (a very mild extension). We already know that $\mathbf{Cs}_b$ is an h1-category, and that 2-homotopy is an equivalence relation (5.3). We also know that concatenation is categorical (associative, with strict identities) and consistent with reversion (2.5.2). We still have to verify the axiom (hc4).

(hc4a) Given a 2-homotopy $\Phi$ and two morphisms $u$, $v$

(1) $\Phi$: $\alpha \simeq_2 \alpha'$: $f \to g$: $X \to Y$, $\qquad u$: $X' \to X$, $\quad v$: $Y \to Y'$,

the following 2-homotopies prove that $v{\circ}\alpha{\circ}u \simeq_2 v{\circ}\alpha'{\circ}u$ and $-\alpha \simeq_2 -\alpha'$:

(2) $P^2v.\Phi.u$: $X' \to P^2Y'$ $\qquad P^2v.\Phi.u$: $v{\circ}\alpha{\circ}u \simeq_2 v{\circ}\alpha'{\circ}u$,

(3) $Pr.\Phi$: $X \to P^2Y$ $\qquad Pr.\Phi$: $-\alpha \simeq_2 -\alpha'$.

Take now a second 2-homotopy $\Psi$: $\beta \simeq_2 \beta'$: $g \to h$: $Y \to Z$. The horizontal faces of $\Phi +_h \Psi$ are 2-homotopic to $\alpha+\beta$, $\alpha'+\beta'$ (5.3), whence $\alpha+\beta \simeq_2 \partial^-P(\Phi +_h \Psi) \simeq_2 \partial^+P(\Phi +_h \Psi) \simeq_2 \alpha'+\beta'$.

(hc4b) To prove that $-\alpha + \alpha \simeq_2 0_g$, take the horizontal pasting of $\Phi = Pr.g^-(\alpha)$ and $\Psi = g^-(\alpha)$. Then $-\alpha + \alpha = \partial^-P.\Phi + \partial^-P.\Psi \simeq_2 \partial^-P(\Phi +_h \Psi) \simeq_2 \partial^+P(\Phi +_h \Psi) = 0_g$.

(hc4c) To prove the weak reduced interchange property

(4) $\varphi\circ f + v\circ\alpha \simeq_2 u\circ\alpha + \varphi\circ g$ (for $\alpha: f \to g: X \to Y$, $\varphi: u \to v: Y \to Z$),

start from the bounded double homotopy $\Phi = P\varphi.\alpha: X \to P^2Z$, $\Phi(i, j; x) = \varphi(i, \alpha(j, x))$, with horizontal faces $\varphi\circ f$, $\varphi\circ g$ and vertical faces $u\circ\alpha$, $v\circ\alpha$; then turn it into a 2-homotopy $\Psi$ realising (4), by "lens conversion" (cf. below, 5.5 a).

**5.5. Lens conversions.** We shall call *lens conversions* the following two procedures, turning *bounded* double homotopies into *bounded* 2-homotopies and vice versa, via connections and (semiregular!) pastings. (The name comes from the topological case [Gr4].)

(a) First, for every bounded double homotopy $\Phi: X \to P^2Y$, with faces as in the left-hand diagram below, there is a bounded 2-homotopy $\Psi: \alpha+\psi \simeq_2 \varphi+\beta: f \to g: X \to Y$

(1)
$$\begin{array}{ccc} k & \xrightarrow{\ \beta\ } & g \\ {\scriptstyle\varphi}\uparrow & \Phi & \uparrow{\scriptstyle\psi} \\ f & \xrightarrow[\ \alpha\ ]{} & h \end{array} \qquad\qquad \begin{array}{ccc} f & \xrightarrow{\ \varphi+\beta\ } & g \\ {\scriptstyle 0}\uparrow & \Psi & \uparrow{\scriptstyle 0} \\ f & \xrightarrow[\ \alpha+\psi\ ]{} & g \end{array}$$

which can be obtained by a horizontal pasting, as in the left-hand diagram below (since concatenation is associative and has strict identities)

(2)
$$\begin{array}{ccccccc} f & \xrightarrow{\ \varphi\ } & k & \xrightarrow{\ \beta\ } & g & \xrightarrow{\ 0\ } & g \\ {\scriptstyle 0}\uparrow & g^+\varphi & \uparrow{\scriptstyle\varphi} & \Phi & \uparrow{\scriptstyle\psi} & g^-\psi & \uparrow{\scriptstyle 0} \\ f & \xrightarrow[\ 0\ ]{} & f & \xrightarrow[\ \alpha\ ]{} & h & \xrightarrow[\ \psi\ ]{} & g \end{array} \qquad\qquad \begin{array}{ccc} k & \xrightarrow{\ \beta\ } & g \\ {\scriptstyle\varphi}\uparrow & g^-\varphi + 0_V(\beta) & \uparrow{\scriptstyle 0} \\ f & \xrightarrow{\ \varphi+\beta\ } & g \\ {\scriptstyle 0}\uparrow & \Psi & \uparrow{\scriptstyle 0} \\ f & \xrightarrow{\ \alpha+\psi\ } & g \\ {\scriptstyle 0}\uparrow & 0_V(\alpha) + g^+\psi & \uparrow{\scriptstyle\psi} \\ f & \xrightarrow[\ \alpha\ ]{} & h \end{array}$$

(b) Conversely, given a bounded 2-homotopy $\Psi: \alpha+\psi \simeq_2 \varphi+\beta$, there exists a bounded double homotopy $\Phi$ with faces $\alpha, \beta, \varphi, \psi$, as in (1); it can be constructed by the right-hand diagram (2) (this second procedure will be used in 5.7).

**5.6. Pseudo-maps and pseudo-homotopies.** We are now able to introduce an *extension* of homotopies which has a concatenation. A *pseudo-map* of simplicial complexes $f = (f_K): X \twoheadrightarrow Y$ will be a family of ordinary maps $f_K: K \to Y$ defined over all finite subspaces $K$ of $X$ (under *no* coherence conditions). The composite with $g = (g_K): Y \twoheadrightarrow Z$ is the family

(1) $(gf)_K = g_{K'}.f_K: K \to Z$ $(K' = f_K(K))$

the identity of $X$ is the family of inclusions $K \to X$. This defines the category Ps**Cs**.

The ordinary category **Cs** will be embedded in Ps**Cs**, preserving objects and identifying a map $f: X \to Y$ with the family of its restrictions $f|_K: K \to Y$ to the finite parts of $X$. The formula (1) reduces to $g.f_K$ when $g$ is a map. The embedding preserves all limits, as it follows easily from the previous formula (or also from the existence of a left adjoint, sending a c-space to the sum of its finite subspaces). But note that Ps**Cs** itself is *not* complete (it lacks equalisers): we shall always use limits of *diagrams of* **Cs**, which are still limits in the extension.

The functor $P: \mathbf{Cs} \to \mathbf{Cs}$ can be extended to Ps**Cs**: given $f = (f_K): X \twoheadrightarrow Y$, take a finite subspace $H \subset PX$; it has a finite "image" $K = \mathrm{Im}(H) = \cup_{a\in H}\, a(\mathbf{Z})$ in $X$, and we let $(Pf)_H$ be the restriction of $P(f_K): K \to PY$ to $H \subset PK$

(2) $(Pf)_H = (P(f_K))|_H: H \to PY, \qquad (Pf)_H(x) = P(f_K)(x);$

a tedious calculation shows it does preserve the composition (1): letting $H \subset PX$ (finite), $K = \mathrm{Im}(H)$, $H' = (Pf)_H(H) = (Pf_K)(H) \subset PY$ and $K' = f_K(K) = \cup_{a\in H}\, f_K a(\mathbf{Z}) = \mathrm{Im}(H') \subset Y$, we have

(3) $(P(gf))_H = (P(gf)_K)|_H = (P(g_{K'}.f_K))|_H = (Pg_{K'}.Pf_K)|_H = (Pg_{K'})|_{H'}.(Pf_K)|_H = (Pg)_{H'}.(Pf_K)|_H =$
$= (Pg.Pf)_H$ .

A *pseudo-homotopy* $\alpha: f \twoheadrightarrow g: X \twoheadrightarrow Y$ will be a pseudo-map $\alpha = (\alpha_K): X \twoheadrightarrow PY$ with $\partial^-\alpha = f$, $\partial^+\alpha = g$; equivalently, it is a family of ordinary homotopies $\alpha_K: f_K \to g_K: K \to Y$, for $K$ finite in $X$. Each of them is bounded, and has been assigned a constant support $\rho(\alpha_K)$. As for ordinary homotopies (3.1.2), the path functor can be extended to pseudo-homotopies, sending $\alpha = (\alpha_K)$ to the family $\hat{P}\alpha = (\beta_H)$ obtained by restriction from $\hat{P}(\alpha_K): P(f_K) \to P(g_K): PK \to PY$

(4) $\beta_H = (\hat{P}\alpha_K)|_H: (Pf)_H \to (Pg)_H: H \to PY \qquad (K = \mathrm{Im}(H) = \cup_{a\in H}\, a(\mathbf{Z})).$

The category Ps**Cs** will be equipped with pseudo-homotopies (represented by the extended path-functor $P$: Ps**Cs** $\to$ Ps**Cs**) and the following structure ($u: X' \to X$, $v: Y \to Y'$)

(a) whisker composition: $(v\circ\alpha\circ u)_K = Pv.\alpha_{uK}.u_K: K \to PY'$ ($K$ finite in $X'$),

(b) trivial homotopies: $(0_f)_K = 0_{f_K}: K \to PY$,

(c) reversion: $(-\alpha)_K = -(\alpha_K): K \to PY$,

(d) concatenation: $(\alpha + \beta)_K = (\alpha_K + \beta_K): K \to PY$,

(e) 2-pseudo-homotopy: $\alpha \simeq_{2p} \beta$ iff ($\alpha_K \simeq_2 \beta_K$, for all $K$ finite in $X$),

where concatenation depends on standard supports and the 2-pseudo-homotopy relation is defined after the similar relation of bounded homotopies (5.3).

All this can be transferred to the pointed case, Ps$\mathbf{Cs}^*$. Then, also *the loop-endofunctor* $\Omega$ *can be extended to pseudo-maps and pseudo-homotopies* (cf. 3.1.5), and is therefore invariant up to pseudo-homotopy: $f \simeq_p g$ implies $\Omega f \simeq_p \Omega g$. This is a crucial fact, showing that the pseudo-homotopy relation still carries a relevant homotopical information.

**5.7. Theorem: The structure of pseudo-homotopies.** Ps**Cs** is a *lax* h4-category with semiregular concatenation, and satisfies (hc1, 3) (as $\mathbf{Cs}_b$, 5.4). The ordinary h-pullback of ordinary maps $P(f, g)$ is a *regular* h-pullback in Ps**Cs**. The same holds for the pointed version, Ps$\mathbf{Cs}^*$.

**Proof.** All axioms are easily deduced from the similar properties of bounded homotopies (5.4). For instance, the weak reduced interchange follows from the same property for bounded homotopies, applied to $\alpha_K: f_K \to g_K$ and $\varphi_L: u_L \to v_L$, where $K$ is finite and $L = \alpha_K(\mathbf{Z}\times K)$ is its (finite) image. An ordinary h-pullback P(f, g) is also a homotopy pullback in Ps**Cs**, because of the characterisation of the latter as an ordinary limit (4.2.2); in fact, PY represents pseudo-homotopies in the category Ps**Cs**, and the embedding preserves all limits (5.6).

Finally, to prove the 2-dimensional property (4.4) of P(f, g) in Ps**Cs**, it suffices to prove the statement for any *finite* simplicial complex W, *in* **Cs**. Start from a 2-homotopy $\Psi$: W → $P^2Z$ representing the 2-homotopy relation of the hypothesis (4.4.1), $\Psi$: $f \circ \varphi + \xi \circ b \simeq_2 \xi \circ a + g \circ \psi$: fpa → gqp: W → Z. Since W is finite, we can modify $\Psi$ by lens conversion (5.5 b), obtaining a double homotopy $\Phi$: W → $P^2Z$ with faces as in 4.4.1 (or 4.5.3, for $a^- = a$, $a^+ = b$). Now, it suffices to apply to $\Phi$ the 2-dimensional property we already have (4.5).

## 6. Intrinsic homotopy groups of simplicial complexes

The interplay of maps, pseudo-maps and pseudo-homotopies allows us to prove the results we are aiming at, for the homotopy groups of pointed simplicial complexes. These are shown to be isomorphic to the homotopy groups of the geometric and McCord realisations of X (6.6).

**6.1. The interplay.** As we have seen in the Structure theorem (5.7), our h-pullbacks produce *ordinary maps* (resp. *pseudo-maps*) from the 1-dimensional property, provided all starting data are ordinary (resp. when some data are in Ps**Cs**); and always *pseudo-homotopies* from the 2-dimensional property. The terms *pseudo-homotopy equivalence* and *pseudo-fibration* will always refer to Ps**Cs**; but the relevant situation is a hybrid one, of *ordinary maps which fall in the previous cases*, within Ps**Cs**.

Let S be a *finite* simplicial complex. Recall that, for maps x: S → X, the various "homotopy relations" ($\simeq_b$, $\simeq_c$, $\simeq$, $\simeq_t$, $\simeq_p$) coincide (3.2); [x] denotes a homotopy class. The resulting functor

(1) [S, –]: **Cs** → **Set**, [S, X] = **Cs**(S, X)/$\simeq$, $f_*$ = [S, f]: [x] ↦ [fx],

is *invariant up to homotopy*: given $\alpha$: f → g, the homotopy $\alpha$x: fx → gx shows that $f_*[x] = g_*[x]$, for all x: S → X. This functor cannot be extended to pseudo-maps f = ($f_K$): X ⇸ Y (as we have no coherence condition on such families, even up to homotopy), but satisfies nevertheless the following invariance properties:

(a) [S, –] is *invariant up to pseudo-homotopy*: for ordinary maps, f $\simeq_p$ g implies $f_* = g_*$.

(b) *an ordinary map* f: X → Y *which is a pseudo homotopy equivalence necessarily induces a bijection* $f_*$: [S, X] → [S, Y].

The first fact is obvious: given x: S → X, we get fx $\simeq_p$ gx, which means that fx $\simeq$ gx (S is finite). For the second, let g = ($g_K$): Y ⇸ X be a quasi-inverse in Ps**Cs**. Given x, x': S → X with fx $\simeq$ fx', it follows that gfx $\simeq_p$ gfx', whence x $\simeq$ x'. Given y: S → Y, take K = y(S) ⊂ Y and x = $g_K$.x: S → X; then fx = f.$g_K$.x $\simeq$ x.

Of course, the same happens in the pointed case, for [S, –]: $\mathbf{Cs}^*$ → $\mathbf{Set}^*$. These results are, again, somewhat parallel to what happens for differential algebras, where a (multiplicative) morphism which becomes a homotopy equivalence in **Dm** (i.e., forgetting the multiplicative structure) induces a multiplicative isomorphism in homology.

The functors of path-components (2.2) are of this type

(2) $\pi_0$: **Cs** $\to$ **Set**, $\quad\quad$ $\pi_0 X = X/\sim = [\{*\}, X]$,

(3) $\pi_0$: **Cs*** $\to$ **Set***, $\quad\quad$ $\pi_0 X = X/\sim = [\mathbf{S}^0, X]$,

where the 0-sphere $\mathbf{S}^0 = \{-1, 1\}$ is discrete, and pointed at 1 (1.3).

**6.2. Homotopy functors.** We have already considered the fundamental groupoid of a simplicial complex X, deduced from the study of paths (2.9). Now, it can be more easily deduced from the structure theorem of bounded homotopies (5.4)

(1) $\Pi$: **Cs** $\to$ **Gpd**, $\quad\quad$ $\Pi X = |PX|/\simeq_2 = (\mathbf{Cs}_b(\{*\}, X)) / \simeq_2$,

since in any lax h4-category **A**, the track 2-category $Ho_2\mathbf{A} = \mathbf{A}/\simeq_2$ is groupoid-enriched, and for each object A, $Ho_2\mathbf{A}(A, -)$ is a functor with values in groupoids.

This functor is *homotopy invariant*, in the sense that a homotopy $\alpha: f \to g: X \to Y$ in **Cs** yields a natural equivalence (a homotopy in **Gpd**, 4.1)

(2) $\Pi\alpha: \Pi f \to \Pi g$, $\quad\quad$ $(\Pi\alpha)x = [\alpha(x)]: f(x) \to g(x)$;

the naturality comes from the weak reduced interchange axiom (in $\mathbf{Cs}_b$), for a path $a: x \to x': \{*\} \to X$ and $\alpha$ itself: $[\alpha(x)] + [ga] = [fa] + [\alpha(x')]$.

Now, let X be a pointed object. $\Omega X \subset PX$ is the object of loops at the base point $*_X$ (2.4.3), and the set $|\Omega X| = |PX|(*_X, *_X)$ is an involutive monoid, under concatenation and reversion.

The (intrinsic) *fundamental group* of X is the group of endomorphisms of $\Pi X$ at the base point. Equivalently, it is the set of path-components of $\Omega X$ (with the operation induced by concatenation)

(3) $\pi_1(X) = (\Pi X)(*_X, *_X) = \pi_0(\Omega X) = |\Omega X|/\sim = [\mathbf{S}^0, \Omega X]$.

(In fact, a path in $\Omega X$, between the loops a, b, amounts to a square $A \in P\Omega X \subset P^2X$ with horizontal faces a, b, subject to the condition that all horizontal paths $A(-, j)$ be loops at $*_X$; but this amounts to saying that the vertical faces of A are trivial, i.e. $a \simeq_2 b$.)

More generally, the (intrinsic) *homotopy functors* of simplicial complexes are defined as

(4) $\pi_n = \pi_0.\Omega^n = [\mathbf{S}^0, \Omega^n(-)]$: **Cs*** $\to$ **Set***,

so that an element $[a] \in \pi_n(X)$ is represented by an n-dimensional net $a: \mathbf{Z}^n \to X$ whose faces are constant at $*_X$, modulo homotopy with fixed faces. All functors $\pi_n$ preserve finite products, since $\pi_0$ and $\Omega$ do (1.5, 2.4); we know they take values in groups for $n \geq 1$.

For $n \geq 2$, we get abelian groups. In fact, the set $|\Omega^2X| \subset |P^2X|$ has two monoid structures, by horizontal and vertical pasting of double paths. Since these two operations have the same identity (the constant double path at the base point) and satisfy the middle four interchange law up to 2-congruence (2.7.6), they induce modulo $\equiv_2$ the same, commutative operation (a well known, easy lemma). As we know that 2-congruence of squares implies homotopy with fixed faces in $P^2X$ (end of 2.8), i.e. pointed homotopy in $\Omega^2X$, the thesis follows.

**6.3. The fibre sequence.** Since **Cs*** is a pointed h1-category, every ordinary map f: X → Y has a fibre sequence (4.3.4)

$$(1) \qquad ... \longrightarrow \Omega Kf \xrightarrow{\Omega x} \Omega X \xrightarrow{\Omega f} \Omega Y \xrightarrow{d} Kf \xrightarrow{x} X \xrightarrow{f} Y$$

connected to the sequence of its iterated h-kernels by a (generally non-commutative) fibre diagram ([Gr2], 6.4).

Since Ps**Cs*** is lax h4, the regularity properties recalled in 4.3 (and proved in [Gr2]) *hold in* Ps**Cs***. Thus, all h-kernel maps are ordinary maps and *pseudo-fibrations*; all comparison maps are ordinary maps and *pseudo-homotopy equivalences*, as well as their $\Omega^n$-images (last remark in 5.6); all comparison squares are *anti-commutative up to pseudo-homotopy* (with respect to the reversion in loop-objects).

If S is finite object, applying to the previous results the functor [S, –], we deduce a sequence of functors still satisfying the properties (a), (b) of 6.1

$$(2) \quad \pi^S_n = [S, \Omega^n(–)]: \mathbf{Cs}^* \to \mathbf{Set}^*, \qquad (\pi^S_n = \pi^S_0\, \Omega^n),$$

and any map f has an exact sequence of pointed mappings (write $\pi_n = \pi^S_n$)

$$(3) \qquad ... \longrightarrow \pi_1 Kf \xrightarrow{\pi_1 x} \pi_1 X \xrightarrow{\pi_1 f} \pi_1 Y \xrightarrow{\pi_0 d} \pi_0 Kf \xrightarrow{\pi_0 x} \pi_0 X \xrightarrow{\pi_0 f} \pi_0 Y.$$

Taking $S = \mathbf{S}^0$, we have the homotopy functors $\pi_n$ (6.2.4) and the *exact homotopy sequence* of the map f. In particular, if f: A → X is the inclusion of a subobject A ≺ X, the relative homotopy "groups" of the pair (X, A) are defined as

$$(4) \quad \pi_n(X, A) = \pi_{n-1}(Kf) = \pi_0(\Omega^{n-1}(Kf)), \qquad (n \geq 1),$$

and (3) becomes the exact homotopy sequence of the pair.

**6.4. A combinatorial van Kampen theorem.** The intrinsic fundamental groupoid ΠX can be calculated by a Seifert-van Kampen-type theorem, similar to the one given by R. Brown for topological spaces ([Br], 6.7); the latter implies the classical result on the fundamental group, but has the advantage of covering more general situations: the circle, to begin with (it does not require a "connected intersection").

If Y is a subobject of the simplicial complex X and R a subset, we shall write $\Pi Y|_R$ the full subgroupoid of ΠY whose objects lie in R∩Y. Say that R is *representative* in Y if it meets all its path-components. Then $\Pi Y|_R$ is a retract of ΠY (choose, for each y∈Y, a Y-path c(y): y → p(y) ∈ R∩Y, and let the functor $p: \Pi Y \to \Pi Y|_R$ send [a]: y → y' to – [c(y)] + [a] + [c(y')]: p(y) → p(y')). Recall also that, given two subobjects U, V ≺ X, the structure of U∪V (resp. U∩V) is !U∪!V (resp. !U∩!V) (1.1); if U and V are subspaces, so are U∪V and U∩V.

*Theorem*. Let X = U∪V be a simplicial complex (we shall say that X *is covered by its subobjects* U, V). If the set R ⊂ |X| is representative in U, V and U∩V, then the following diagram of groupoids (induced by the inclusions) is a pushout

$$(1) \qquad \begin{array}{ccc} \Pi(U\cap V)|_R & \longrightarrow & \Pi U|_R \\ \downarrow & & \downarrow \\ \Pi V|_R & \longrightarrow & \Pi X|_R \end{array}$$

The proof is quite similar to the one in [Br], simplified by the fact that here any n-tuple path $a: \mathbf{Z}^n \to X$ has a *standard* decomposition over the elementary cubes $\xi$ of its support in $\mathbf{Z}^n$: for each of them, $a(\xi)$ is linked in $U$ or in $V$; when it is so in both, then it is also linked in $U\cap V$. (In the topological case one constructs a non-standard decomposition, invoking the Lebesgue covering theorem).

**6.5. The combinatorial circles.** It is now easy to deduce that the circles $C_k$ (1.3) have fundamental group $\mathbf{Z}$. (Alternatively, one can recall that their geometric realisation is the circle, and apply the isomorphism through the edge-path groupoid $\Pi C_k \cong \mathcal{E} C_k \cong \Pi\mathbf{S}^1$, 2.10).

If $k > 3$, $X = C_k$ is covered by two *subspaces* $U = \{[0], [1], [2]\}$, $V = \{[2],... [k–1], [0]\}$, which are contractible (being isomorphic to integral intervals) and whose intersection is the discrete object on two points. Taking $R = |U\cap V| = \{[0], [2]\}$, $\Pi U|_R$ is the groupoid with two objects and two reciprocal isomorphisms connecting them, $\bullet \rightleftarrows \bullet$ ; $\Pi V|_R$ is the same, while $\Pi(U\cap V)|_R$ is the discrete groupoid at the two selected points. Thus, we have obtained the same three groupoids which one finds for the *topological* circle, with the same embeddings. Their pushout $\Pi X|_R$ is isomorphic to the fundamental groupoid of the latter (at two points), and the fundamental group of $X$ is the free group on one generator, represented by the obvious loop $[0], [1],... [k]$.

*For* $C_3$ *we need to work with subobjects*. The subobject $U$ has a tolerance structure on the same set $|C_3| = \{[0], [1], [2]\}$, with $[0] ! [1] ! [2]$, while $V = \{[2], [0]\} \subset C_3$ has the induced (chaotic) structure. Then, $U$ and $V$ are contractible and $U\cap V$ is discrete; one concludes as before.

The van Kampen pushout proves also that all maps $p_k: C_{k+1} \to C_k$ induce isomorphism on the fundamental groupoids and $\pi_1$. But our combinatorial circles are *not* homotopy equivalent: if $k < k'$, the image of any map $C_k \to C_{k'}$ is a proper, connected subspace, isomorphic to an integral interval, hence contractible.

While $\pi_0 X = [\mathbf{S}^0, X]$ is "homotopy representable", this is not true for the higher homotopy functors: **Cs*** *has no standard sphere in positive dimension*. Working, for simplicity, in dimension 1, the problem reduces to the fact that a given loop cannot "represent longer ones". More precisely, suppose by absurd that there be a simplicial complex $S$ and a natural isomorphism $[S, X] = \pi_1 X = |\Omega X|/\!\sim$ (the brackets denote equivalence classes for $\simeq_t$). Then, the identity $idS$ determines a special class in $\pi_1 S$, represented by some loop $a_0: \mathbf{Z} \to S$, and an arbitrary class $[a] \in \pi_1 X$ corresponds to a homotopy class $[f]$, where $f: S \to X$ is a map; but $f = f\circ idS = [S, f](idS)$, and by naturality $[a] = (\pi_1 f)[a_0] = [fa_0]$. Thus, for any $X$, any class $[a] \in \pi_1 X$ has a representative whose support is contained in $\rho(a_0)$. This is absurd, since each representative of the generator of $\pi_1(C_k)$ is surjective.

However, the whole system of maps $p_k: C_{k+1} \twoheadrightarrow C_k$ surrogates a "standard circle", in the sense that $\pi_1 X$ can be expressed as a colimit (i.e., a direct limit)

(1) $\pi_1 X = \mathrm{colim}_k\, [C_k, X]$, $\qquad [p_k, X]: [C_k, X] \to [C_{k+1}, X]$,

the isomorphism being induced by the standard loops $q_k: \mathbf{Z} \twoheadrightarrow C_k = \mathbf{Z}/\!\equiv_k$, $q_k(i) = [(0\vee i)\wedge k]$

(2) $[C_k, X] \to \pi_1 X$, $[a: C_k \to X] \mapsto [aq_k: \mathbf{Z} \to X]$.

This suggests that it might be useful to embed **Cs** and **Cs*** in their categories of pro-objects (cf. [MaS, CP]), and realise there a "true" circle; and similarly a standard interval.

**6.6. Theorem: The simplicial comparison.** Let $X$ be a pointed simplicial complex, $\mathcal{R}X$ its geometric realisation. Then, there is a natural isomorphism $\pi_n(X) \to \pi_n(\mathcal{R}X)$ between the combinatorial and topological homotopy groups.

**Proof.** Following the description of [Sp], $\mathcal{R}X$ is the set of all mappings $\lambda: X \to [0, 1]$ with linked support $\mathrm{supp}(\lambda)$, such that $\Sigma_x \lambda(x) = 1$. $X$ is embedded in $\mathcal{R}X$, identifying $x \in X$ with its characteristic function. A point of $\mathcal{R}X$ can be viewed as a convex combination $\lambda = \Sigma_i \lambda_i x_i$ of a linked family of $X$ and each (non-empty) linked subset $\xi$ having $p+1$ points spans a *simplex*

(1) $\Delta(\xi) = \{\lambda \in \mathcal{R}X \mid \mathrm{supp}(\lambda) \subset \xi\}$, $\Delta(\xi) \supset \Delta°(\xi) = \{\lambda \in \mathcal{R}X \mid \mathrm{supp}(\lambda) = \xi\}$.

All $\Delta(\xi)$ are equipped with the euclidean topology (via a bijective correspondence with the standard simplex $\Delta^p$, determined by any linear order of $\xi$), and $\mathcal{R}X$ with the direct limit topology defined by such subsets. Each $\Delta(\xi)$ is closed in $\mathcal{R}X$, whereas the *open simplex* $\Delta°(\xi)$ is open *in* $\Delta(\xi)$; in particular, $\{x\} = \Delta(\{x\}) = \Delta°(\{x\})$. $\mathcal{R}X$ is the disjoint union of all its open simplices.

For each $x \in X$, the *star*

(2) $\mathrm{st}(x) = \{\lambda \in \mathcal{R}X \mid x \in \mathrm{supp}(\lambda)\} = \cup_{\xi \ni x} \Delta°(\xi)$, (disjoint union),

is an open neighbourhood of $x$ in $\mathcal{R}X$, and their family is an open covering. Then

(3) $\mathrm{st}(x) \cap X = \{x\}$, $\mathrm{st}(y) \cap \Delta(\xi) \neq \emptyset$ iff $y \in \xi$.

A subset $\xi \subset X$ is linked iff $\cap_{x\in\xi} \mathrm{st}(x) \neq \emptyset$, and then

(4) $\mathrm{st}(\xi) = \cap_{x\in\xi} \mathrm{st}(x) = \{\lambda \in \mathcal{R}X \mid \xi \subset \mathrm{supp}(\lambda)\} = \cup_{\eta\supset\xi} \Delta°(\eta) \supset \Delta°(\xi)$;

note that $\mathrm{st}(\xi)$ and $\Delta°(\xi)$ do not meet $\xi$, as soon as the latter has at least two points. But $\mathrm{st}(\xi)$ is "star-shaped with respect to $\Delta°(\xi)$", i.e. each convex combination $t.\lambda + (1–t).\mu$ $(0 \leq t \leq 1)$ with $\lambda \in \mathrm{st}(\xi)$ and $\mu \in \Delta°(\xi)$ is in $\mathrm{st}(\xi)$ (because its support is linked, $\mathrm{supp}(\lambda)$ or $\xi$).

The image of a combinatorial mapping $a: \mathbf{2}^n \to X$ is a linked subset $\xi$ of $X$, contained in the convex space $\Delta(\xi) \subset \mathcal{R}X$, and we can consider the *multi-affine* extension $\hat{a}: [0, 1]^n \to \mathcal{R}X$ of $a$ (separately affine in each variable), here exemplified for $n = 2$

(5) $\hat{a}(s, t) = (1–s)(1–t).a(0, 0) + s(1–t).a(1, 0) + (1–s)t.a(0, 1) + st.a(1, 1)$.

Every net $a: \mathbf{Z}^n \to X$ has a unique extension to a continuous mapping $\hat{a}: \mathbf{R}^n \to \mathcal{R}X$, affine in each variable on all elementary cubes $\Pi_r [i_r, i_r+1]$ with vertices in $\mathbf{Z}^n$ (and eventually constant, at the left and the right, in each variable). Moreover, if certain faces of $a$ are trivial (constant at the base point), so are the corresponding faces of $\hat{a}$. If $\zeta = \Pi\zeta_r \subset \mathbf{Z}^n$ is linked, we shall write $\Delta(\zeta)$ the convex envelope of $\zeta$ in $\mathbf{R}^n$ and $\Delta°(\zeta)$ the interior of $\Delta(\zeta)$ in the affine subspace which it generates; then, $\hat{a}(\Delta(\zeta)) \subset \Delta(a(\zeta))$ and $\hat{a}(\Delta°(\zeta)) \subset \Delta°(a(\zeta))$. This defines a natural homomorphism

(6) $\Phi_n: \pi_n(X) \to \pi_n(\mathcal{R}X)$, $\Phi_n[a] = [\hat{a}]$.

To prove its surjectivity, take a (Moore) n-dimensional topological path $f: [0, \rho]^n \to \mathcal{R}X$, with trivial faces. By the Lebesgue covering theorem (and by continuity at the faces of the cube), one can find a partition $0 = t_0 < t_1 < ... < t_{k+1} = \rho$ of the interval $[0, \rho]$ and a family $x_i \in X$ such that

(7) $f(\Pi_r [t_{i_r-1}, t_{i_r+1}]) \subset st(x_i)$, $\quad i = (i_1,... i_n) \in \mathbf{Z}^n, \ 1 \leq i_1,... i_n \leq k$,

$x_i = *_X$, $\quad$ when some $i_k$ is 1 or k,

so that $f(\Pi_r [t_{i_r}, t_{i_r+1}]) \subset st(x(\zeta))$, for $\zeta = \Pi_r [i_r, i_r+1]$.

Prolonging this family $x_i$, we have a net $x: \mathbf{Z}^n \to X$ (with support $[1, k]^n$ and trivial faces) since an elementary cube $\zeta_i = \{i+j \mid j\in \mathbf{2}^n\}$ is sent to a subset $\xi_i = \{x_{i+j} \mid j\in \mathbf{2}^n\}$, which is linked because the stars of its elements meet

(8) $f(t_{i_1},... t_{i_n}) \in f(\Pi_r [t_{i_r+j_r-1}, t_{i_r+j_r+1}]) \subset st(x_{i+j})$ $\quad (i \in [1, k-1]^n, \ j \in \mathbf{2}^n)$.

Modifying the parametrisation of f, we can assume that $f: [0, k+1]^n \to \mathcal{R}X$ and $t_j = j$. The continuous mapping $\hat{x}: \mathbf{R}^n \to \mathcal{R}X$ is n-homotopic to f, by a homotopy affine in the deformation parameter; in fact, if $t\in [0, k+1]^n$, then $t\in \Delta°(\zeta)$ for some $\zeta = \Pi\zeta_r$ linked in $\mathbf{Z}^n \cap [0, k+1]^n$; thus

(9) $\hat{x}(t) \in \Delta°(x(\zeta))$, $\quad f(t) \in st(x(\zeta))$,

(the second fact is analogous to the statement which follows (7)) and we know that $st(x(\zeta))$ is star-shaped with respect to its convex subset $\Delta°(x(\zeta))$.

As to injectivity, take two nets $a, b: \mathbf{Z}^n \to X$ with support $[0, \rho]^n$ and trivial faces, together with a map $f: [0, \rho]^{n+1} \to \mathcal{R}X$, with last faces $\hat{a}$, $\hat{b}$ and all other faces trivial. The latter determines a partition $0 = t_0 < t_1 < ... < t_k = \rho$ such that $f(\Pi_r [t_{i_r-1}, t_{i_r}]) \subset st(x_i)$ and a net $x: \mathbf{Z}^{n+1} \to X$ as above, all whose faces except the last ones are trivial. This is a combinatorial n-homotopy between its last faces $\partial^{\kappa}_{n+1}(x)$, and it is sufficient to show that these are n-homotopic to a, b, respectively.

Assume, for simplicity, that $t_j - t_{j-1} \leq 1$ (whence $k \geq \rho$). Thus, the mapping $d: \mathbf{Z} \to \mathbf{Z}$

(10) $d(j) = j \quad (j \leq 0)$, $\quad d(j) = [t_j] \quad (0 \leq j \leq k)$, $\quad d(j) = j + \rho - k \quad (j \geq k)$,

is a delay (2.6); a is n-homotopic to the net $a\mathbf{d} = a(d\times...\times d)$, which is immediately n-homotopic (3.1) to $\partial^-_{n+1}(x) = (x_{i0})$ $(i = (i_1,...i_n))$, since an elementary cube $\zeta_i = \{i+j \mid j\in \mathbf{2}^n\}$ is sent by both nets into a common linked subset

(11) $a\mathbf{d}(i+j) \in \{a(\mathbf{d}(i)+j')) \mid j'\in \mathbf{2}^n\}$, $\quad x_{i+j,0} \in \{a(\mathbf{d}(i)+j')) \mid j'\in \mathbf{2}^n\}$.

The first fact is obvious; the second follows from (3): $st(x_{i+j,0})$ contains the point $f(t_{i_1},... t_{i_n}, 0)$, by (8), which also belongs to the convex envelope of the previous linked subset

(12) $f(t_{i_1},... t_{i_n}, 0) = \hat{a}(t_{i_1},... t_{i_n}) \in \Delta\{a(\mathbf{d}(i)+j')) \mid j'\in \mathbf{2}^n\}$,

by definition of $\hat{a}$ and because $d(i_r) \leq t_{i_r} \leq d(i_r)+1$.

## 7. Combinatorial homotopy of metric spaces

We study the fundamental group $\pi_1^\varepsilon(X)$ of a metric space, at resolution $\varepsilon$. Computations are performed by means of the combinatorial van Kampen theorem (6.4) and telescopic homotopies (Section 3); or, in some cases, can be reduced to the topological $\pi_1$ of the *spot dilation* (7.4-5).

**7.1. Some calculations.** Let X be a pointed metric space. Recall that, for $0 \leq \varepsilon \leq \infty$, we have considered the tolerance structure $t_\varepsilon X$, with x!x' iff $d(x, x') \leq \varepsilon$, and the homotopy "groups" at resolution $\varepsilon$, forming a direct system of groups (or pointed sets, for n = 0)

(1) $\pi_n^\varepsilon(X) = \pi_n(t_\varepsilon X)$, $\qquad \pi_n^\varepsilon(X) \to \pi_n^\eta(X) \qquad (0 \leq \varepsilon \leq \eta \leq \infty)$;

for $\varepsilon = 0$ (resp. $\infty$) $t_\varepsilon X$ is totally disconnected (resp. chaotic, hence contractible) and all $\pi_n^\varepsilon(X)$ are trivial for $n > 0$. Let us take $0 < \varepsilon < \infty$.

We are now able to compute $\pi_1^\varepsilon(X)$ for the metric space 1.8.1, $X = T \setminus (A \cup B \cup C)$, pointed at the origin (say), a closed region of the real plane (with the $l_\infty$-metric), endowed with the $t_\varepsilon$-structure

(2) 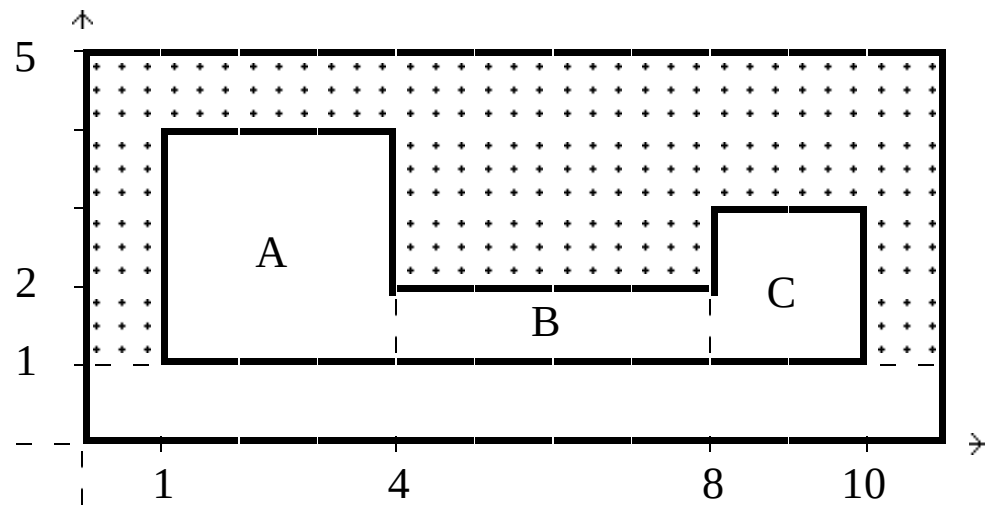

$T = [0, 11] \times [0, 5]$

$A = ]1, 4[ \times ]1, 4[$
$B = [4, 8] \times ]1, 2[$
$C = ]8, 10[ \times ]1, 3[$.

To apply the combinatorial van Kampen theorem (6.4), we shall use a variable subspace $U \subset X$ and a fixed $V = \{(x, y) \in X \mid y > 1\}$ (the dotted region above); the latter is $t_\varepsilon$-telescopic with respect to the line y = 5, with constant jump $s_i = \varepsilon$ (3.6), hence contractible.

(a) If $0 < \varepsilon < 1$, the argument proceeds as in 6.5. X is covered by its subspaces $U = \{(x, y) \in X \mid y < 2\}$ and V; U is contractible (use again a constant jump $s_i = \varepsilon$), while $U \cap V = ([0, 1] \times ]1, 2[) \cup ([10, 11] \times ]1, 2[)$ has two connected components, which are contractible. Taking as a representative subset $R = \{(0, 3/2), (10, 3/2)\}$, one deduces that the fundamental group $\pi_1^\varepsilon(X)$ is the free group on one generator, represented by a loop around $A \cup B \cup C$.

(b) Let $1 \leq \varepsilon < 2$. Cover X with $U = \{(x, y) \in X \mid y < 3\}$ and V; then U, the dotted subspace at the left hand below, is (connected and) $t_\varepsilon$-telescopic with respect to the first axis (with jumps $s_i = 1$)

(3) 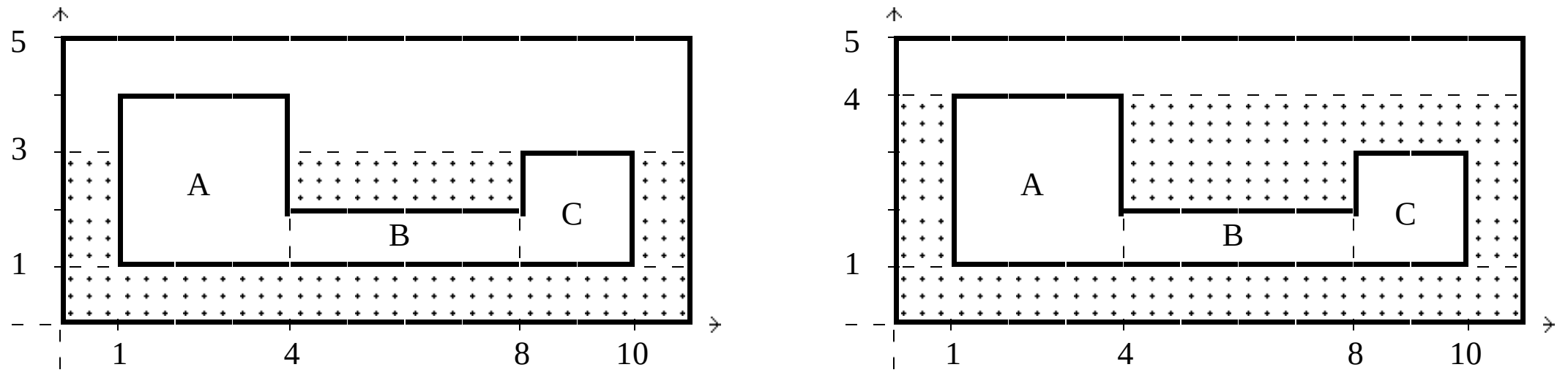

Thus U, V are contractible, while $U \cap V$ has three connected components, contractible. Choosing a representative point in each of them, e.g. $R = \{(0, 2), (4, 2), (10, 2)\}$, one falls in the same computation of the topological fundamental groupoid of the "figure 8"; thus, $\pi_1^\varepsilon(X)$ is the free group on two generators, represented by two loops, one around A and the other around C.

(c) Let $2 \le \varepsilon < 3$. Cover X with $U = \{(x, y) \in X \mid y < 4\}$ and V; now U (the dotted subspace at the right hand, above) is $t_\varepsilon$-telescopic with respect to the axis $y = 0$, with vertical jumps $s_1 = 1$, $s_2 = 2$, $s_3 = 1$ (to jump over $B \cup C$). Thus, U, V are contractible and $U \cap V$ has two connected components, contractible. $\pi_1^\varepsilon(X)$ is the free group on one generator, represented by a loop around A.

(d) For $\varepsilon \ge 3$, we already showed that X is contractible (3.6).

(e) The same arguments work on the finite model $[X]_\rho = |X| \cap \rho\mathbf{Z}^2$ ($\rho = k^{-1}$, for a positive integer k) considered in 1.8, provided that $\varepsilon \ge \rho$ (the only modification is to use a constant jump $s_i = \rho$ to show that V is contractible, as well as for U in the case $\rho \le \varepsilon < 1$); for $\varepsilon < \rho$, we have a totally disconnected object.

**7.2. Critical values.** Such results lead us to consider the *variation* of the system of homotopy groups $\pi_n^\varepsilon(X) \to \pi_n^\eta(X)$ $(0 \le \varepsilon \le \eta \le \infty)$, for a pointed metric space X, through its *critical values* (cf. Deheuvels [De], Milnor [Mi]). Here we can just give some hints.

Say that $\varepsilon$ is a *left regular* (resp. *right regular*, *regular*) value for this system (or for X, in dimension n) if the system itself is constant on a left (resp. right, bilateral) neighbourhood of $\varepsilon$ in the extended real interval $[0, \infty]$. In the contrary, $\varepsilon$ is a *left critical* (resp. *right critical*, *critical*) value; and a *bilateral critical* value if it is both left and right critical.

Thus, the pointed metric space X considered in 7.1 has in dimension 1 a right critical value at 0 and left critical values at 1, 2, 3 (and no else). The metric subspace $X = T \setminus (A \cup B) \subset \mathbf{R}^2$ represented below (pointed at the origin) has in dimension 1

(1) 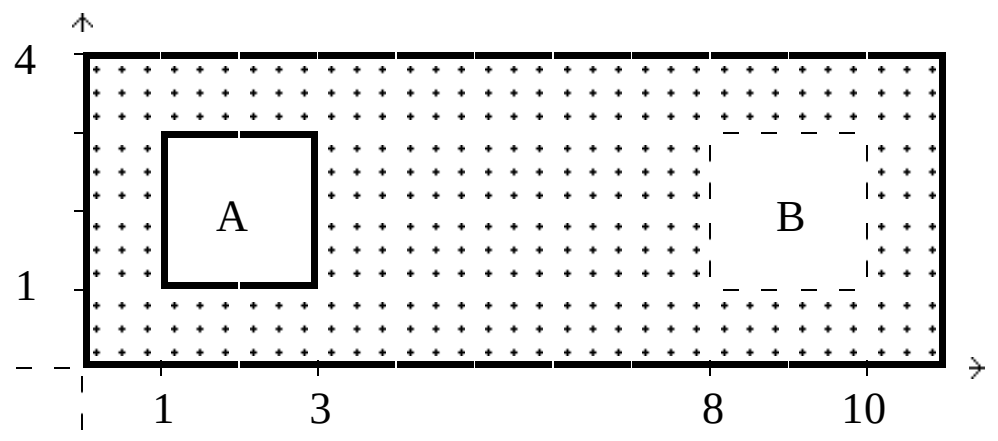

$T = [0, 11] \times [0, 4]$
$A = ]1, 3[ \times ]1, 3[$
$B = [8, 10] \times [1, 3]$

a right critical value at 0 and a bilateral critical value at 2. In fact, its fundamental $\varepsilon$-group can be calculated by van Kampen and telescopic homotopies, as above

(2) $\pi_1^\varepsilon(X) \;\cong\; \mathbf{Z} * \mathbf{Z}$ $(0 < \varepsilon < 2)$, $\qquad \mathbf{Z}$ $(\varepsilon = 2)$, $\qquad \{*\}$ $(\varepsilon > 2)$.

For $0 < \varepsilon < 2$, the argument proceeds as in 7.1b. For $\varepsilon = 2$, X is covered by its subspaces U, V determined by the inequalities $x \le 10$ and $x \ge 8$; both are contractible (they are telescopically homotopic, with suitable jumps, with respect to the lines $x = 0$ and $x = 11$, respectively), and $U \cap V$ has two connected components, which are contractible. For $\varepsilon > 2$, X is contractible, being telescopically homotopic with respect to the line $8 - \rho$ $(0 < 2\rho < (\varepsilon - 2) \wedge 2)$, with characteristic sequence

(3) $a_{-4} = 1, \quad a_{-3} = 3, \quad a_{-2} = 5, \quad a_{-1} = 7, \quad a_0 = 8 - \rho, \quad a_1 = 10 + \rho.$

**7.3. Open structures.** Let us come back to a general metric space X. For $\varepsilon > 0$, the "open" tolerance structure $t_\varepsilon^- X$ defined by $d(x, x') < \varepsilon$ would seem at first to be of the same interest as $t_\varepsilon X$, since each family determines the other, in the complete lattice of tolerance structures on the set $|X|$

(1) $t_\varepsilon^- X = \cup_{\eta<\varepsilon} t_\eta X,$ $\qquad$ $t_\varepsilon X = \cap_{\eta>\varepsilon} t_\eta^- X,$

yet the first family $(t_\varepsilon X)_\varepsilon$ yields finer results in homotopy. In fact, its homotopy groups determine the others

(2) $\pi_n(t_\varepsilon^- X) = \mathrm{colim}_{\eta<\varepsilon}\, \pi_n(t_\eta X),$

as a trivial consequence of the finiteness of paths: a net $a: \mathbf{Z}^n \to t_\varepsilon^- X$ of the new structure is also a map $\mathbf{Z}^n \to t_\eta X$, with $\eta = \mathrm{diam}(a) = \max_\xi(\mathrm{diam}(a(\xi))) < \varepsilon$, where $\xi$ varies in the finite set of elementary cubes of the support of a.

On the other hand, there cannot be any general formula giving $\pi_n^\varepsilon(X)$ as a limit or colimit of the other family $\pi_n(t_\varepsilon^- X)$, since in the previous example (7.2.1) the family $\pi_1(t_\varepsilon^- X)$, calculated directly or deduced from the colimit (2), takes only two values: $\mathbf{Z}*\mathbf{Z}$ $(0 < \varepsilon \leq 2)$ and $\{*\}$ (otherwise).

**7.4. Spot dilations.** We end with a construction which realises simplicial complexes on regions X of $\mathbf{R}^n$ as *topological subspaces* of $\mathbf{R}^n$ itself and often allows for an easy calculation of the fundamental group: a "dilation operator", as used in mathematical morphology (cf. [He, FM]), which expands X to a union of discs. Coming back to previous remarks in 1.9, the efficiency of this procedure for the c-spaces considered here should be compared with the hugeness of their geometric realisation, infinite-dimensional in the continuous cases; on the other hand, of course, the geometric realisation is general and – theoretically – works for all homotopy groups, while the present comparison is limited (necessarily) to particular metric c-spaces and (presently) to $\pi_0$, $\pi_1$.

Let $\varepsilon > 0$ and let X be a pointed subset of a normed real vector space E, with metric d. The *open-spot dilation* $D_\varepsilon^-(X)$ and the *closed-spot dilation* $D_\varepsilon(X)$ of X in E will be the following (pointed) metric subspaces, obtained as a union of d-discs *of radius* $\varepsilon/2$

(1) $D_\varepsilon^-(X) = D_\varepsilon^-(X, d) = \{x' \in E \mid d(x, x') < \varepsilon/2 \text{ for some } x \in X\} \supset X,$

(2) $D_\varepsilon(X) = D_\varepsilon(X, d) = \{x' \in E \mid d(x, x') \leq \varepsilon/2 \text{ for some } x \in X\} \supset X.$

At present, we are only able to give a general result for $\pi_1$ (and $\pi_0$) and the *first* operator, relative to *both* tolerance structures $t_\varepsilon^- X$ and $t_\varepsilon X$. Say that X is $t_\varepsilon^-$-*closed* (resp. $t_\varepsilon$-*closed*) in E if $D_\varepsilon^-(X)$ contains the convex envelope of all linked subsets of $t_\varepsilon^- X$ (resp. $t_\varepsilon X$), so that there is a continuous mapping

(3) $f_\varepsilon^-: \mathfrak{S}(t_\varepsilon^- X) \to D_\varepsilon^-(X)$ $\qquad$ (resp. $f_\varepsilon: \mathfrak{S}(t_\varepsilon X) \to D_\varepsilon^-(X)$),

extending the identity of X and affine on each simplex $\Delta(\xi)$ of the domain. We prove below that there is then a canonical isomorphism $\pi_1(t_\varepsilon^- X) \cong \pi_1(D_\varepsilon^- X)$ (resp. $\pi_1(t_\varepsilon^- X) \cong \pi_1(t_\varepsilon X) \cong \pi_1(D_\varepsilon^- X)$), induced by this map via the geometric comparison (6.6). Of course, $t_\varepsilon$-closure implies $t_\varepsilon^-$-closure.

The applications are straightforward. Take $E_p = (\mathbf{R}^2, d_p)$ with the $l_p$-norm $(1 \leq p \leq \infty)$ and its continuous regions X of 7.1.2 or 7.2.1. Then, X is $t_\varepsilon^-$-closed in $E_p$, for all p and all $\varepsilon$, which allows one to calculate quite easily the fundamental group $\pi_1(t_\varepsilon^-(X, d_p))$ as $\pi_1(D_\varepsilon^-(X, d_p))$. X is also $t_\varepsilon$-closed in $E_p$, for $p < \infty$ and all $\varepsilon$ *different from* 1, 2, 3 (different from 2 for 7.2.1). Moreover $\mathbf{Z}^2$ is $t_\varepsilon^-$-closed in $E_1 = (\mathbf{R}^2, d_1)$, for $1 < \varepsilon \leq 2$

(4) $t^-_\varepsilon(\mathbf{Z}^2, d_1) \; = \; t_1(\mathbf{Z}^2, d_1) \; = \; t_1(\mathbf{Z}^2, d_p)$ $\qquad (1 \leq p < \infty, \; 1 < \varepsilon \leq 2),$

$D^-_2(\mathbf{Z}^2) \; = \; \mathbf{R}^2 \setminus \{(m + 1/2, n + 1/2) \mid (m, n) \in \mathbf{Z}^2\}$ $\qquad$ in $(\mathbf{R}^2, d_1)$,

showing that $\pi_1(t_1(\mathbf{Z}^2, d_1)) = \pi_1(D^-_2(\mathbf{Z}^2))$ is a free group of countable rank (while we know that $\mathbf{Z}^2 = t_1(\mathbf{Z}^2, d_\infty)$ is contractible).

Coming back to the general case, it would be interesting to have a similar result for the closed-spot dilation, i.e. an isomorphism $\pi_1(t_\varepsilon X) \cong \pi_1(D_\varepsilon X)$ under sufficiently general hypotheses, satisfied by the regions considered in 7.1 or 7.2 (where this coincidence appears to hold for the $l_\infty$-norm, *without exceptions at critical values*). Here it is not sufficient to assume that $D_\varepsilon(X)$ contain the convex envelope of all linked subsets of $t_\varepsilon X$. In fact, for $\varepsilon = 1$, this condition holds for $X = \{0\} \cup ]1, 2] \subset \mathbf{R}$; but $t_1 X$ has two path components, while the metric space $D_1 X = [-1/2, 5/2]$ is connected; it is easy to deduce a 2-dimensional counterexample, with different $\pi_1$. But the previous condition might be sufficient for a *compact* X, which would include cases as 7.1.2.

**7.5. Theorem: The open-spot comparison.** Let X be a pointed metric subspace of a normed vector space E and $\varepsilon > 0$.

(a) If X is $t^-_\varepsilon$-closed (7.4), there is a canonical isomorphism $\Psi^-_1 \colon \pi_1(t^-_\varepsilon X) \to \pi_1(D^-_\varepsilon X)$ between the combinatorial and topological homotopy groups, which is the composite of the geometric realisation isomorphism $\Phi_1$ (6.6.6) with the *iso*morphism induced by the "affine" map $f^-_\varepsilon$ (7.4.3)

(1) $\qquad \Psi^-_1 \; = \; \pi_1(f^-_\varepsilon).\Phi_1 \; = \; (\pi_1(t^-_\varepsilon X) \; \to \; \pi_1(\mathfrak{Z}(t^-_\varepsilon X)) \; \to \; \pi_1(D^-_\varepsilon X)).$

(b) If X is also $t_\varepsilon$-closed, there is a canonical isomorphism $\Psi_1 \colon \pi_1(t_\varepsilon X) \to \pi_1(D^-_\varepsilon X)$ consisting of the lower row of the following commutative diagram of isos

$$
(2) \qquad
\begin{array}{ccccc}
\pi_1(t^-_\varepsilon X) & \xrightarrow{\;\Phi_1\;} & \pi_1(\mathfrak{Z}(t^-_\varepsilon X)) & \xrightarrow{\;\pi_1(f^-_\varepsilon)\;} & \pi_1(D^-_\varepsilon X) \\
\downarrow & & \downarrow & & \| \\
\pi_1(t_\varepsilon X) & \xrightarrow{\;\Phi_1\;} & \pi_1(\mathfrak{Z}(t_\varepsilon X)) & \xrightarrow{\;\pi_1(f_\varepsilon)\;} & \pi_1(D^-_\varepsilon X)
\end{array}
$$

where the vertical arrows are induced by the subobject inclusion $t^-_\varepsilon X \prec t_\varepsilon X$.

In (a) and (b), similar conclusions hold for $\pi_0$.

**Proof.** (a) Since $f^-_\varepsilon \colon \mathfrak{Z}(t^-_\varepsilon X)) \to D^-_\varepsilon X$ is affine on simplices, the composed homomorphism $\Psi^-_n = \pi_n(f^-_\varepsilon).\Phi_n$ acts much as $\Phi_n$ in 6.6: every n-dimensional net $a \colon \mathbf{Z}^n \to t^-_\varepsilon(X)$ has a unique extension to a map $\hat{a} \colon \mathbf{R}^n \to D^-_\varepsilon(X)$, affine in each variable on all elementary cubes $\prod_r [i_r, i_r{+}1]$ with vertices in $\mathbf{Z}^n$ (and eventually constant, the left and the right, in each variable). We prove directly that $\Psi^-_1$ is iso, by an argument similar in part to the one of 6.6; its extension for $n > 1$ is not evident, while the case $n = 0$ is obvious.

For surjectivity, take a Moore loop $f \colon [0, \rho] \to D^-_\varepsilon(X)$. By the Lebesgue covering theorem (and by continuity at the end points), we can find two finite sequences $t_i \in [0, \rho]$, $x_i \in X$ such that

(3) $\quad 0 \; = \; t_0 \; < \; t_1 \; < \; \dots \; < \; t_k \; = \; \rho, \qquad f([t_{i-1}, t_i]) \; \subset \; D(x_i, \varepsilon/2), \qquad x_0 \; = \; x_k \; = \; *_X,$

where $D(x, \lambda)$ denotes the open disc of E of centre x and radius $\lambda$. The sequence $x_i$ (prolonged) is a combinatorial loop $x \colon \mathbf{Z} \to t^-_\varepsilon X$, since

(4) $d(x_i, x_{i+1}) \leq d(x_i, f(t_i)) + d(f(t_i), x_{i+1}) < \varepsilon$,

and the topological loop $\hat{x}$ is equivalent to the original f in $D_\varepsilon^-(X)$, by a homotopy affine in the time parameter, since the (convex) disc $D(x_i, \varepsilon/2) \subset D_\varepsilon^-(X)$ contains $f([t_{i-1}, t_i])$ together with the points $(x_{i-1} + x_i)/2$, $(x_i + x_{i+1})/2$ of $\hat{x}$.

For injectivity, take two combinatorial loops $a, b: \mathbf{Z} \to t_\varepsilon^-(X)$ with support $[0, \rho]$ and a topological 2-homotopy $f: [0, \rho]^2 \to D_\varepsilon^-(X)$, with horizontal faces $\hat{a}$, $\hat{b}$. We can similarly find two finite families $t_i \in [0, \rho]$, $x_{ij} \in X$ $(0 \leq i, j \leq k)$ such that

(5) $0 = t_0 < t_1 < ... < t_k = \rho$, $\qquad f([t_{i-1}, t_i]\times[t_{j-1}, t_j]) \subset D(x_{ij}, \varepsilon/2)$,

$x_{0j} = x_{kj} = *_X$.

By refining this partition, we can assume that:

(i) $t_i$ assumes all *semi-integral* values $m + 1/2$ $(m\in \mathbf{Z})$ of $[0, \rho]$ (whence $t_i - t_{i-1} \leq 1$),

(ii) if $t_i = m + 1/2$ $(m\in \mathbf{Z})$, then $x_{i-1,j} = x_{ij}$ (for all j).

Extending the family $(x_{ij})$ to $x: \mathbf{Z}^2 \to t_\varepsilon^-(X)$, we obtain a combinatorial 2-homotopy between its horizontal faces $(x_{i0})$ and $(x_{ik})$. The mapping $\delta: \mathbf{Z} \to \mathbf{Z}$

(6) $\delta(i) = i \quad (i \leq 0)$, $\qquad \delta(i) = [t_i + 1/2] \quad (0 \leq i \leq k)$, $\qquad \delta(i) = i + \rho - k \quad (i \geq k)$,

is a delay (because $t_i - t_{i-1} \leq 1$); thus a is 2-homotopic to the path $a\delta$, and it is sufficient to prove that the latter is immediately 2-homotopic to $(x_{i0})$, i.e. that each four-tuple $a\delta(i)$, $a\delta(i+1)$, $x_{i0}$, $x_{i+1,0}$ is linked $(0 \leq i < k)$. We shall show that there is a point $y\in E$ sufficiently near:

(7) $d(x_{i0}, y) < \varepsilon/2$, $\qquad d(x_{i+1,0}, y) < \varepsilon/2$,

$d(a\delta(i), y) = d(f(\delta(i), 0), y) < \varepsilon/2$, $\qquad d(a\delta(i+1), y) = d(f(\delta(i+1), 0), y) < \varepsilon/2$ .

By assumption (i), it cannot happen that $t_i < m+1/2 < t_{i+1}$ $(m\in \mathbf{Z})$. Thus

(8) $m - 1/2 \leq t_i < t_{i+1} \leq m + 1/2$, $\qquad \delta(i) = m, \quad \delta(i+1) = m$ or $m+1$.

First, suppose that $t_{i+1} < m+1/2$, so that $\delta(i+1) = m = \delta(i)$, and take $y = f(t_i, 0)$. Then, the first and second inequality of (7) follow from (5). The third (and fourth) follows from the fact that $f(t, 0) = \hat{a}(t)$ is affine on each interval $[n, n +1]$, with $d(a(n), a(n+1)) < \varepsilon$, and the numbers $t_{i+1}$ and $\delta(i) = m$ fall in a common half of such intervals, either $[m - 1/2, m]$ or $[m, m + 1/2]$. Finally, suppose that $t_{i+1} = m+1/2$, so that $\delta(i+1) = m+1$ and $x_{i0} = x_{i+1,0}$ (by ii); take now $y = f(t_{i+1}, 0)$. The second (and first) inequality of (7) follows from (5); the third and fourth follow again from $f(t, 0)$ being affine on the interval $[m, m+1] = [\delta(i), \delta(i+1)]$, of which $t_{i+1}$ is the middle point.

(b) If X is also $t_\varepsilon$-closed, we have the commutative diagram (2). We already know that the two geometric comparisons $\Phi_1$ are isos, as well as $\pi_1(f_\varepsilon^-)$. The same argument as in (a) proves that the composed morphism $\Psi_1$ of the lower row is also so. Therefore, all arrows of (2) are isos.